\documentclass[10pt,openany,leqno]{amsart}
\usepackage[backref=page,colorlinks=true]{hyperref}
\usepackage[capitalize]{cleveref}
\usepackage{amsmath,amsthm,amsfonts,amssymb,amscd,url}
\usepackage{graphics}
\numberwithin{equation}{section}
\input{xy} \xyoption{all}
\usepackage[utf8]{inputenc}

\usepackage{enumerate,xcolor}

\usepackage{epsfig}
\input{xy} \xyoption{all}

\usepackage[scr=boondoxo,scrscaled=1.05]{mathalfa}
\makeindex
\begin{document}
\newcommand{\bfemph}[1]{\emph{\bf{#1}}}
\def\Card{\mathrm{Card\, }}
\def\E{\mathrm{Symp}^\infty_{el}}
\def\cT{\mathrm{Symp}^\infty_{T}}
\def\Leb{\mathrm{Leb\,}}
\def\exp{\mathrm{exp}}
\def\S{\mathbb S}
\def \mC{{i\rho}}
\def\u{\underline }
\def\inv{^{-1}}
\def\Max{\text{max}}
\def\cal{\mathcal}
\def\Graph{\mathrm{Graph\, }}
\def \sfw{{\mathsf w}}
\def \sfh{{\mathsf h}}
\def\R{\mathbb R}
\def\N{\mathbb N}
\def\Z{\mathbb Z}
\def\D{\mathbb D}
\def\A{\mathbb A}
\def\B{\mathbb B}
\def\C{\mathbb C}
\def\O{\mathbb O}
\def\cC{\mathcal C}
\def\T{\mathbb T}
\def\a{{\underline a}}
\def\b{{\underline b}}
\def\c{{\underline c}}
\def\Log{\mathrm{log}}
\def\loc{\mathrm{loc}}
\def\inta{\mathrm{int }}
\def \Diff{\mathrm{Diff}}
\def \Emb{\mathrm{Emb}}
\def\det{\mathrm{det}}
\def\Re{\mathrm{Re}}
\def\lip{\mathrm{Lip}}
\def\leb{\mathrm{Leb}}
\def\dom{\mathrm{Dom}}
\def\diam{\mathrm{diam}\:}
\def\supp{\mathrm{supp}\:}
\newcommand{\ovfork}{{\overline{\pitchfork}}}
\newcommand{\ovforki}{{\overline{\pitchfork}_{I}}}
\newcommand{\Tfork}{{\cap\!\!\!\!^\mathrm{T}}}
\newcommand{\whforki}{{\widehat{\pitchfork}_{I}}}
\newcommand{\marginal}[1]{\marginpar{{\scriptsize {#1}}}}
\def \np{{\color{red} 4}}
\def \npm{{\color{red} 3}}
\def \npmm{{\color{red} 2}}
\def \det{\mathrm{det}\, }
\def\cO{{\mathcal O}}
\def\cU{{\mathcal U}}
\def\cB{{\mathcal B}}
\def\cD{{\mathcal D}}
\def\cG{{\mathcal G}}
\def\cP{{\mathcal P}}
\def\cF{{\mathscr F}}
\def\cH{{\mathcal H}}
\def\cV{{\mathcal V}}
\def\cW{{\mathcal W}}
\def\cX{{\mathcal X}}
\def\cY{{\mathcal Y}}
\def\cZ{{\mathcal Z}}
\def\cE{{\mathcal E}}

\def\cN{{\mathcal N}}

\def\cM{{\mathcal M}}
\def\qand{\quad \text{and} \quad}

\def\sA{{\mathfrak A}}
\def\sB{{\mathfrak B}}
\def\sC{{\mathfrak C}}
\def\sE{{\mathfrak E}}
\def\sG{{\mathfrak G}}
\def\sI{{\mathfrak I}}
\def\sJ{{\mathfrak J}}
\def\sQ{{\mathfrak Q}}
\def\sP{{\mathfrak P}}
\def\sR{{\mathfrak R}}
\def\sS{{\mathfrak S}}
\def \sT{{\mathfrak T}}
\def\sU{{\mathfrak U}}
\def\sV{{\mathfrak V}}
\def\sX{{\mathfrak X}}
\def\sY{{\mathfrak Y}}
\def\sZ{{\mathfrak Z}}

\def\Symp{\mathrm{ Symp}}

\def\sa{{\mathfrak a}}
\def\sb{{\mathfrak b}}
\def\sc{{\mathfrak c}}
\def\sd{{\mathfrak d}}
\def\se{{\mathsf e}}
\def\sf{t}
\def\sg{{\mathfrak g}}
\def\sh{{\mathfrak h}}
\def\si{o}
\def\sm{{\mathfrak m}}
\def\sn{{\mathfrak n}}
\def\sq{{\mathfrak q}}
\def\Lip{\mathrm{Lip}}

\def\so{{\diamond}}
\def\sp{{\mathfrak p}}
\def\sl{{\mathfrak {l}}}
\def\sr{{\mathfrak {r}}}
\def\ss{{\mathfrak {s}}}
\def\st{{\mathfrak {t}}}
\def\su{{\mathfrak {u}}}
\def\sv{{\mathfrak {v}}}
\def\spp{{\mathfrak {p}}}
\def\sw{{\mathfrak {w}}}
\def\arr{\overleftarrow}
\def\avv{\overrightarrow}
\def\epsilon{\varepsilon}
\def\eps{\varepsilon}
\theoremstyle{plain}
\def\Cr#1{\overset{\tilde Y}{#1}}

\newtheorem{thm}{\bf Theorem}[section]
\newtheorem{theorem}[thm]{\bf Theorem}
\newtheorem*{theorem*}{\bf Theorem}

\newtheorem*{conjecture*}{\bf Conjecture}

\newtheorem{conj}[thm]{\bf Conjecture}
\newtheorem{claim}[thm]{\bf Claim}
\newtheorem{assumption}[thm]{\bf Assumption}

\newtheorem{problem}[thm]{\bf Problem}
\newtheorem{question}[thm]{\bf Question}
\newtheorem{proposition}[thm]{\bf Proposition}
\newtheorem{corollary}[thm]{\bf Corollary}
\newtheorem{lemma}[thm]{\bf Lemma}

\newtheorem{sublemma}[thm]{\bf Sublemma}
\newtheorem*{Takens prbm}{Takens' Last Problem} 
\newtheorem{remark}[thm]{\bf Remark}
\newtheorem{fact}[thm]{\bf Fact}
\newtheorem{exem}[thm]{\bf Example}
\newtheorem{definition}[thm]{\bf Definition}
\newtheorem*{definition*}{\bf Definition}
\newtheorem{defi}[thm]{\bf Definition}

\newtheorem{theo}{\bf Theorem} 
\newtheorem{coro}[theo] {\bf Corollary}
\renewcommand{\thetheo}{\Alph{theo}}
\renewcommand{\thecoro}{\Alph{theo}}

\newtheorem{example}[thm]{\bf Example}
\renewcommand*{\backref}[1]{}
\renewcommand*{\backrefalt}[4]{\quad \tiny 
  \ifcase #1 (\textbf{NOT CITED.})%
  \or    (Cited on page~#2.)%
  \else   (Cited on pages~#2.)%
  \fi}

\title{On Kolmogorov-typical properties of symplectic dynamics}

\author{Pierre Berger} 
 
\author{Dmitry Turaev}

\thanks{Pierre Berger is  partially supported by the ERC project 818737 Emergence of wild differentiable dynamical systems
 and by the French National Research Agency under the project DynAtrois (ANR-24-CE40-1163)}
\thanks{Dmitry Turaev is supported by the Leverhulme Trust grant RPG-2021-072}
\address{Pierre Berger\\
Institut de Mathématiques de Jussieu-Paris Rive Gauche\\
CNRS, UMR 7586\\
Sorbonne université, Université Paris Cité\\
4 place Jussieu, 75252 Paris Cedex 05}
\address{Dmitry Turaev\\ Imperial College, London, SW72AZ UK}

\begin{abstract}We propose a general framework, within which we prove that several properties, such as the fast growth of the number of periodic points, the universality, and the high emergence, hold true for every parameter value for a generic finite-parameter family of symplectic diffeomorphisms displaying an elliptic point.
\end{abstract}
\maketitle
\tableofcontents

\section{Introduction}
The study of symplectic dynamics began with the Hamiltonian systems. Most of physical observations of such systems display elliptic (i.e., stable in the linear approximation) periodic orbits.
Non-degenerate elliptic orbits are accumulated by a fibration by Lagrangian invariant tori, called KAM-tori (where KAM stands for Kolmogorov-Arnold-Moser \cite{Kol54, Arn63, Mo62}). The union of these tori is a set of positive, and relatively large, Lebesgue measure and corresponds to quasiperiodic motions. The complement to this set of regular dynamics is poorly understood. In this work, we show from several view points that a neighborhood of an elliptic orbit \emph{typically} contains, in the complement to the set of KAM-tori, invariant regions
with dynamics of the highest complexity and richness.

To make precise statements, we employ the notion of {\em Kolmogorov typicality} in the spirit of  \cite{Kol57,IL99,Ber19}. 
Let $(M,\omega)$ denote a compact symplectic manifold of dimension $2n$ and let $\Symp^\infty(M)$ denote the Fr\'echet space of $C^\infty$-smooth symplectomorphisms of $M$.  
Let $\E(M)$ be the open subset of $\Symp^\infty(M)$ of symplectomorphisms displaying a non-degenerate elliptic periodic point (the non-degeneracy means here the standard twist conditions
for the Birkhoff normal form, see \cref{def_non_deg}). A property~$(\cP)$ is \bfemph{Kolmogorov typical} in $\E(M)$ if, given \emph{any} compact manifold $\sA$ of 
parameters $a$ and a {\em $C^\infty$-generic} family $(f_a)_{a\in \sA}$ of symplectomorphisms $f_a\in \E(M)$, the property $(\cP)$ holds for the maps $f_a$ at \emph{every}  parameter value $a$.

This definition of the Kolmogorov typicality is actually stronger than usual, as we precise and discuss in \cref{state of art}: it does \emph{not} depend on $\dim \sA\ge 0$.
 
\subsection{Characterizations of the dynamical complexity} 
We derive the following Theorems \ref{Main explicit per},  \ref{Main explicit entropy}, \ref{Main explicit universality} and \ref{Main explicit emergence} from a more abstract Theorem \ref{theo Main} stated in section \ref{localmain}.
\begin{theo}\label{Main explicit per} Kolmogorov-typically, a symplectomorphism $f\in \E(M)$ has its number of $n$-periodic points $Per_n(f)$ growing super exponentially:
\[ \limsup_{k\to \infty}\frac{ \log \# Per_k(f) }{\log k} =+ \infty\; .\]
\end{theo}
\begin{theo}\label{Main explicit entropy} There is a dense set of smooth families $(f_a)_{a \in \sA}$  of symplectomorphisms $f_a\in \E(M)$ such that $f_a$ has positive metric entropy for every $a\in \sA$:
\[ \int \liminf_{k\to \infty} \frac1k  \log \| Df^k _a  \| d\leb >0 \; .\]
\end{theo}
\begin{example}We recall that by Duarte's Theorem \cite{Du94},  the celebrated Chirikov standard map $S_a: (x,p)\mapsto (x + p + a \sin x, p + a \sin x)$  displays non-degenerate elliptic points for an open and dense set of large parameters $a$.  Hence, Theorem \ref{Main explicit entropy} implies that a $C^\infty$-small perturbation of the family $(S_a)_{a\in \R} $ has a positive metric entropy for a non-empty open set of parameters $a$.  \end{example}

The next application regards the notion of universality. A map $f\in \Symp^\infty(M)$ is \bfemph{universal} \cite{Btalk, BD03,Tu03},
if the set of its renormalized iterations is $C^\infty$-dense in the space of smooth symplectic diffeomorphisms from the closed $2n$-dimensional unit ball $\B$ onto its image in $\R^{2n}$, where 
$2n=\dim M$. A \bfemph{renormalized iteration} of $f$ is a map $F$ such that there exist $k\ge 1$ and a $C^\infty$-smooth, conformally-symplectic embedding 
$h: B\cup F(B)\hookrightarrow M$ such that $h\circ F|_B =  f^k\circ  h|_B$. The conformal symplecticity of the coordinate transformation $h$ means it takes the standard symplectic form $\omega$ in $\R^{2n}$ to 
$\lambda \omega$ for some constant scalar $\lambda\neq 0$; thus the renormalized iteration $F$ preserves the standard symplectic form. Since $h(B)$ can have an arbitrarily small size and be situated anywhere in $M$, the set of renormalized iterations represents the entirety of local dynamics of the given map $f$ at arbitrarily fine spatial scales. By definition, if a map is universal, then its renormalized iterations $C^\infty$-approximate all bounded symplectic dynamics in $\R^{2n}$. 
Thus, the following implies that the ultimate richness of the fine-scale dynamics is Kolmogorov-typical among symplectomorphisms with elliptic orbits:
\begin{theo}\label{Main explicit universality} Kolmogorov-typically, a symplectomorphism from $\E(M)$  is universal. 
\end{theo}

The \bfemph{emergence}  \cite{Be16} describes how large is the distribution of possible 
statistics generated by a dynamical system, i.e., how large is the ergodic decomposition $\se_*\Leb$ of 
the Lebesgue measure on $M$. If a symplectic map is ergodic with respect to the Lebesgue measure, then 
$\se_*\Leb=\delta_{\leb}$ is atomic. If the map is the identity, then $ \se_*\Leb=\int_M  \delta_{x}d\leb$. 
In \cite{BB21},  it was shown that a locally generic surface diffeomorphism needs a super polynomial number of  probability measures to $\epsilon$-approximate its ergodic decomposition (the quantization dimension is infinite).
Then, one estimates the size of the ergodic decomposition by its {\em order}.
The \bfemph{emergence order is maximal} if a sample of approximately 
$\sim \exp\,  {\epsilon^{-\dim M}}$ probability measures is needed to $\epsilon$-approximate
$\se_*\Leb$ (i.e., to describe the statistical behavior of the orbits with 
accuracy $\epsilon$).  In this case, the order of $\se_*\Leb$ is approximately as large as the order of the space of all probability measures on $M$; the notion of maximal order for emergence is made precise in \cref{section def emergence}.
According to the following theorem, the ultimate complexity of the statistical behavior is Kolmogorov-typical for symplectomorphisms with elliptic orbits:
\begin{theo}\label{Main explicit emergence} Kolmogorov-typically, a symplectomorphism from $\E(M)$ has emergence of maximal order.\end{theo}

\subsection{Kolmogorov typicality of generic localizable properties}\label{localmain}
Our main theorem employs the concept of a \bfemph{localizable property} that goes back to \cite{Ber19}.
We discuss properties of $C^\infty$-symplectomorphisms of a compact symplectic manifold $(M,\omega)$, but the construction makes also sense for other classes of diffeomorphisms.
\begin{definition}\label{localizable equiv prop}
A property $(\cP)$ is localizable if for any nonempty open set $U\subset  M$ there exists a close to the identity, compactly supported $g\in \Symp^\infty (U)$ such that  the property $(\cP)$
holds for any $f\in \Symp^\infty(M)$ such that
$f^N|_U=g$ for some $N\geq 1$.
\end{definition}
The following implies that having a positive metric entropy is a localizable property:
\begin{proposition}\label{prop entropy}
The property of having a positive Lyapunov exponent 
over a subset of positive Lebesgue measure is localizable in 
$\Symp^\infty(M)$. 
\end{proposition} 
\begin{proof}In \cite{BT19}, we  showed that a compactly supported $C^\infty$-perturbation $g_0$ of the identity map of the disk have non-zero Lyapunov exponent $\lambda=\liminf_{k\to \infty} \frac1k  \log \| Dg_0^k(x)\|$
for the set of points $x$ of positive Lebesgue measure. Taking the product of $n$-copies of $g_0$ gives the proposition.
\end{proof}

A property $(\cP)$ on $\Symp^\infty(M)$ is \bfemph{openly localizable}, if every map   in a neighborhood in $\Symp^\infty(M)$ of the set of maps $f$ described in \cref{localizable equiv prop} satisfies property $(\cP)$. 
 For example, we have immediately:
\begin{proposition}\label{prop open per}
Let $(c_k)_k\in \N^\N$ be any sequence of integers. The property of having more than $c_k$ of $k$-periodic points for a certain $k$ is openly localizable in  $\Symp^\infty(M)$.\end{proposition}
The following is non trivial and is proved in \cref{Universality proof}:
\begin{proposition}\label{prop open renorm}
Let  $\cU\subset \Symp^\infty(\B\!\hookrightarrow\! \R^{2n})$ be a nonempty open subset of symplectic embeddings. Then the property of having a renormalized iteration belonging to $\cU$ is  openly localizable in  $\Symp^\infty(M)$.\end{proposition}
We notice that the conjunction  of two  openly localizable properties is openly localizable. Indeed it suffices to extract from $U$  two disjoint open subsets $U_1$ and $U_2$ and apply the definition to both sets.  A property $(\cP)$ is \bfemph{generically localizable}, if it is implied by the conjunction of a countable number of openly localizable properties $\{(\cP_i): {i\ge 0}\}$:
\[\bigwedge_i (\cP_i)\Rightarrow (\cP)\; .\]
For instance, one has from \cref{prop open per} that having a fast growth of the number of periodic points is a generically localizable property (see also  \cite[Prop. 2.2]{Ber19}). Likewise \cref{prop open renorm} implies that the universality is a generically localizable property.

Using the KAM theorem, we will prove in \cref{emergence proof}:
 \begin{proposition}\label{example of localizable prop3}
Having emergence of maximal order is a generically localizable property in $\Symp^\infty(M)$.
\end{proposition}

The above propositions immediately imply Theorems \ref{Main explicit per}, \ref{Main explicit entropy}, \ref{Main explicit universality} and \ref{Main explicit emergence} according to Theorem \ref{theo Main}, which we state now. 
 
Given a compact smooth manifold $\sA$, we denote as $\E(M)_\sA$ the space of $C^\infty$-families  $f_\sA= (f_a)_{a\in \sA}$ of maps $f_a\in \E(M)$ that we endow with the topology induced by $C^\infty (A\times M, M)$ via the inclusion 
\[(f_a)_{a \in \sA}\in \E(M)_\sA\hookrightarrow
  [(a,x)\mapsto f_a(x)]\in C^\infty(\sA\times M, M).\]
Then $\E(M)_\sA$  is a Fr\'echet space and so Baire. A  \emph{generic subset} of $\E(M)_\sA$ is a countable intersection of open-dense subsets. 
\begin{theo}[Main]\label{theo Main}Let $(\cP)$ be a property on the space of symplectomorphisms of  $(M,\omega)$.
\begin{enumerate}
\item If $(\cP)$ is localizable, then for every compact manifold  $\sA$, the set of families  
$(f_a)_{a\in \sA}$ such that $(\cP)$ holds for the maps $f_a$ at every $a\in \sA$ is dense in $\E(M)_\sA$. 
\item If $(\cP)$ is openly locallizable, then for every compact manifold  $\sA$, the set of  families $(f_a)_{a\in \sA}$ such that
$(\cP)$ holds for $f_a$ at every $a\in \sA$ comprises an open and dense subset of $\E(M)_\sA$. 
\item If $(\cP)$ is generically locallizable, then for every compact manifold  $\sA$, the set of  families  
$(f_a)_{a\in \sA}$ such that $(\cP)$ holds for $f_a$ at every $a\in \sA$ comprises a generic subset of $\E(M)_\sA$, 
i.e., $(\cP)$ is  Kolmogorov typical in $\E(M)$. 
\end{enumerate} 
\end{theo}

We derive this theorem from Theorem \ref{main1} below. For now, let us mention that Theorem \ref{theo Main} is in line with the following general
\begin{conjecture*}[\cite{Ber19}]\label{Conjecture}
There are many classes of dynamical systems for which any generic localizable property is Kolmogorov typical.
\end{conjecture*}
Let us also formally deduce Theorems \ref{Main explicit per}, \ref{Main explicit entropy}, \ref{Main explicit universality} and \ref{Main explicit emergence} from  Theorem \ref{theo Main}. 
\begin{proof}[Proof of Theorem \ref{Main explicit per}] Let $\cU$ be an open subset of $M$ and $N\ge 1$. Define 
$c^{(m)}_k= \exp \circ \exp( k+m)$.  By \cref{prop open renorm},  for any given $m\geq 0$, the property $(\cP_{m})$ of having more than $c^{(m)}_k$ periodic points of period $k$ for some $k\ge 1$ is openly localizable. Taking the intersection over 
all $m\ge 0$, we obtain that the property of displaying a super-exponential growth of the number of periodic points is generically localizable. Then the third item of Theorem \ref{theo Main} implies that this property is Kolmogorov typical in $\E(M)$. 
\end{proof}
\begin{proof}[Proof of Theorem \ref{Main explicit entropy}] It is a direct consequence of \cref{prop entropy}  and   item (1) of Theorem~\ref{theo Main}.\end{proof}
\begin{proof}[Proof of Theorem \ref{Main explicit universality}] It follows from \cref{prop open renorm} and item (3) of Theorem~\ref{theo Main}.\end{proof}
\begin{proof}[Proof of Theorem \ref{Main explicit emergence}] It is immediately implied by \cref{example of localizable prop3}  and   item (3) of Theorem~\ref{theo Main}.\end{proof}

We prove Theorem \ref{theo Main} by the analysis of local bifurcations near an elliptic periodic point. 
In \cref{accum KAM}, we recall the definition of a non-degenerate elliptic point and remind the basic fact that these points are accumulated by \emph{persistent KAM tori}. Such a torus is a Lagrangian $n$-torus, invariant with respect to an iterate $f^q$ of a  map, and this iteration acts as a Diophantine rotation on the torus. Moreover, the dynamics in the normal bundle of such a torus is non-degenerate. This enables bringing $f^q$ in a neighborhood of the KAM torus, by a symplectic coordinate transformation, to the following Birkhoff normal form for any integer $D\geq 1$:
\[ f^q (\theta, r) = (\theta+ \alpha + \nabla Q(r), r)+ O(r^{D+1}),\]
where $(\theta,r)\in \T^n \times \B^n$, the $n$-vector $\alpha$ is Diophantine, $Q$ is a degree $(D+1)$ polynomial such that $\nabla Q(0)=0$ and $\det \nabla^2 Q(0)\neq 0$.  According to the KAM theory (see \cref{KAM classique}), when we unfold $f$ in a family $(f_a)_{a\in \sA}$, the invariant torus with the frequency vector $\alpha$ persists, and the 
coordinate transformation that brings $f^q$ near the torus to the Birkhoff normal form depends smoothly on $a$, see \cref{Birkhoff}. 

Thus, by a $C^{D}$-small perturbation of the family $(f_a)_{a\in \sA}$ within $\E(M)_\sA$, for a smooth family of symplectic  coordinates we have
\[ f_a^q(\theta, r) =  (\theta+ \alpha + Q_a(r), r)\]
near the KAM torus, for an open set of values of $a$ and, importantly, $\alpha$ staying constant as $a$ changes. By an additional small perturbation, we can make $\alpha$ to be a rational vector, so all points in the torus $r=0$ become periodic and, in a neighborhood of each of these periodic points, an iterate of the dynamics gets locally conjugate to a linear parabolic map.  Hence, it can be perturbed to a rational rotation. Note that  a rational rotation displays  an invariant open subset $U$  on which an iterate of  the dynamics is the identity. This line of reasoning will lead us to obtain the following result in \cref{proof main1}:
\begin{theorem}\label{main1}In any neighborhood of any family $(f_a)_{a \in \sA}\in  \E(M)_\sA$,  there exists a family $(\tilde f_a)_{a \in \sA}$ for which 
 there is an open covering $(\sA_i)_{1\le i\le I}$  of $\sA$ and for every $1\le i\le I$, there are  a nonempty open subset $U_i$ and an integer $N_i\ge 1$ such that 
 $f^{N_i}_a|_{U_i}=id$ for every $a\in \sA_i$.
\end{theorem}
From this, we deduce Theorem \ref{theo Main}.
\begin{proof}[Proof of Theorem \ref{theo Main}]
Up to reducing $N_i$, we can assume that there exists a point $x$ in $U_i$ whose minimal period is  $N_i$.  This is also true for every point close to $x$. Hence, up to reducing $U_i$, we can assume that $f_a^k(U_i)$ is disjoint from $U_i$ for every $1\le k<N_i$.
Also,  up to reducing a second time the $U_i$, we can assume that $f_a^k(U_j)$ is disjoint from $f_a^\ell (U_i)$ for all $i\neq j$ and  $k\le N_j$, $\ell \le N_i$. This implies, as property $(\cP)$ is localizable, that there exists a perturbation of $(\tilde f_a)_a$ localized at $\sqcup_i U_i$ such that $f_a$  satisfies~$(\cP)$ for every $a\in \sA$. This gives $(1)$. Then (2) is obvious. 

Let us prove $(3)$.  Let $(\cP_i)_{i\ge 0} $ be the localizable, open properties whose conjunction implies 
$(\cP)$. By $(2)$, for every $i\ge 0$,  there is an open and dense set of families in $\E(M)_\sA$ which satisfies $(\cP_i)$.   Thus all the properties are satisfied at a countable intersection of open and dense set of families,  which is indeed a  generic set of families in $\E(M)_\sA$. 
\end{proof}

\subsection{State of the art}\label{state of art} 

The enhanced complexity of chaotic dynamics due to the violation of hyperbolicity was discovered by Newhouse in a breakthrough paper \cite{Ne74},
see also \cite{Ne79}. He showed that  for an open set $\cN$ of $C^r$-diffeomorphisms, $2\le r\le \infty$, a generic diffeomorphism from $\cN$ displays  infinitely many periodic sinks and the closure of the set of these sinks contains hyperbolic Cantor sets of positive topological entropy. This result was transferred to surface symplectomorphisms (with replacing sinks by elliptic islands) by Duarte  \cite{Du94,Du08}; similar results for the inseparable coexistence of infinitely many saddles of different indices, sources, and also non-trivial attractors of various types were obtained in \cite{GTS93a, Bu97, GTS97, Co98,  BD99,GTS08,GTS09,Tu15, GT17, BGG23, BFP23}. 
Furthermore, such phenomena have been shown to appear at every parameter value in generic families \cite{BE15,BCP22}. 
 
The mechanism behind the Newhouse phenomenon and its generalizations is the existence of robust homoclinic tangencies. In contrast, the robust mechanism used for the main result of this work is the persistence of KAM tori. These can be perturbed to produce ``periodic spots''  (domains filled by periodic orbits only), which in turn can be perturbed to produce homoclinic tangencies and, in fact, dynamics of arbitrary complexity. Such scheme enables to show that  various complex phenomena are typical in a very strong sense.\\
 
\paragraph{\bf  Typicality} 
During his plenary talk at the ICM 1956, Kolmogorov proposed to distinguish between the properties which hold  \emph{typicaly} among 
Hamiltonian dynamical systems, and those which are satisfied  only for a negligible set of systems. Yet, he warned that this question is ``extremely tricky''. He considered as \emph{typical} a property which happens for a system $f_0$, as well as  almost every map $f_a$ close enough to $f_0$ 
in  an analytic unfolding  $(f_a)_{a\in (-\epsilon, \epsilon)}$.  He also conjectured that this notion should have an interesting finitely differentiable counterpart. These notions of typicality have inspired several studies  \cite{IL99,KH07,KH072,Ber19}, to quote only a few. 
 Another notion of typicality, is the topological genericity.  A property is said to be topologically generic on a Fr\' echet space, if it holds at a countable intersection of open and dense subsets. Note that the notion of Kolmogorov typicality we deal with is much stronger: the families can be parametrized by \emph{any} compact manifold $\sA$ and we ask the property to hold true at \emph{every} parameter.  In particular  if a property is Kolmogorov typical, then it is topologically generic (take $\sA$ equal to a point).\\ 

\paragraph{\bf Emergence}  One of the cornerstones of statistical physics is the ``Gibbs postiulate'', which prefers the Liouville measure to all other invariant measures of a Hamiltonian system. The ergodicity of the Liouville measure would provide a support for this postulate, but the ergodic hypothesis is known to be false since the seminal work of Kolmogorov on near-integrable dynamics \cite{Kol54}: a positive measure invariant set filled by KAM-tori supports uncountably many ergodic components.
The ubiquity of elliptic islands in Hamiltonian dynamics -- also away from the near-integrability \cite{Du94}, as well as in low- and high-energy multi-particle systems \cite{Rink, RKT24} -- makes the robust ergodicity of the Liouville measure implausible without the partial hyperbolicity constraint \cite{PS97,ACW21}.

Yet, the mere presence of a fibration by KAM tori does not intuitively appear as a strong obstruction for the understanding of the statistical behavior of a system. Indeed, only a polynomial number $\approx \epsilon^{-n}$  of statistics suffices to describe the dynamics in such  integrable region with the precision~$\epsilon$. In order to understand whether the statistics of generic Hamiltonian dynamics is more complicated than that, we employ the notion of ``emergence'' that has been introduced in \cite{Be16} to quantify the complexity of the statistical behavior of the ``wild'' systems (i.e., those displaying infinitely many attractors or, in the symplectic setting, ergodic components). 
    
The emergence is high if the number of  probability measure needed to describe the system up to the precision $\epsilon$ is super polynomial (Sup-P). Following the celebrated Cobham's thesis \cite{Cob65}, the statistical description of any system with the Sup-P high emergence is not algorithmically feasible.

The local topological genericity of high emergence has been shown in  \cite{BB21} for surface diffeomorphisms and surface symplectomormisms. On the other hand, examples in a very rigid class of systems  and with high emergence have been constructed in \cite{BB23} among surface polynomial automorphism of $\R^2$ and in \cite{Ber22} among analytic surface symplectomorphisms. See also \cite{De25} for other surfaces.
 
Theorem \ref{Main explicit emergence} is the first result on typicality of high emergence in the sense of Kolmogorov. It proves the following conjecture in  the category of symplectomorphisms: 
\begin{conj}[\cite{Be16}]\label{the conj}
Typical dynamics display a super polynomial emergence, in many senses and many contexts.
\end{conj}
A stronger notion is that of the local emergence, see \cite{Ber19}. Having a \emph{high local emergence} means that the local dimension $\limsup_{\epsilon\to 0}\frac{ \se_* \leb (B(\mu,\epsilon))}{ \log \epsilon}$   of $ \se_* \leb$ is not finite on a set of positive measure. It can be shown that it is a localizable property; an important open problem is:
\begin{problem} Does high local emergence is a  localizable generic property?
\end{problem}

Another view on  statistical complexity (Turing uncomputability) is studied e.g. in \cite{RY24}.\\

\paragraph{\bf Universality} The concept of the universal dynamics (as Bonatti put it, ``dynamics containing the whole Universe'' \cite{Btalk})
has been developed since
\cite{BD03} and \cite{Tu03} (in resp. $C^1$ and $C^\infty$ contexts). A dynamical system is universal in a class of maps if the set of its renormalized iterations approximates all maps of the class. Typically we speak about maps (diffeomorphisms or endomorphisms) of a $d$-dimensional ball, for different values of the dimension $d$ and different regularities, or with additional structures (symplectomorphisms, volume-preserving maps,
reversible maps).  In \cite{AST16}, universal iterated function systems were also constructed. The universality by itself is a manifestation of the utmost complexity of the dynamics; it can be utilized for applied purposes, e.g. in machine learning \cite{J+20,V+22}.

For several cases, the genericity of universal dynamics has been established e.g. in \cite{GTS93,GTS08,GeT10,Tu15,AST16,AST20,GT17,BFP23,Mi25}. The main tool of this theory is the bifurcation analysis of homoclinic tangencies and heterodimensional cycles, see \cite{Tu10}. An important ingredient is the theorem from \cite{Tu06,Tu15} that universal diffeomorphisms exist arbitrarily close, in $C^\infty$, to the identity map of a ball (a stronger version of this result has been obtained in \cite{BGH}; the analysis of the perturbations of the identity map goes back to the celebrated papers by Ruelle, Takens and Newhouse \cite{RT,NRT}). This theorem reduces the problem of the search of the universal dynamics to the search of periodic spots; this approach is implemented in the present paper in the symplectic setting. Earlier, periodic spots near elliptic points of surface symplectomorphisms were found in \cite{GeT10}, and near KAM-tori of reversible maps in \cite{GT17} (in any dimension).\\

\paragraph{\bf Positive metric entropy}The positive metric entropy conjecture states that a ``typical'' surface symplectomorphism has positive metric entropy. Here ``typical'' may be interpreted as occurring in some specific natural examples, like the Chirikov standard maps, as conjectured by Sinai \cite[P. 144]{Si94}, or over a large subset of dynamical systems (i.e., dense or generic).  In \cite{BT19}, we proved Herman's conjecture \cite{He98} by showing that any surface symplectomorphism displaying an elliptic periodic point can be $C^\infty$-approximated by one displaying positive metric  entropy. Our theorem was proved by showing that the dynamics constructed by surgery from a uniformly hyperbolic system (similar to a Przytycki example \cite{Pr82}), appears after a perturbation, and a subsequent renormalization, of the identity map of a disc. These examples can be realized among  analytic surface symplectomorphisms \cite{Ber23} and for perturbations of an integrable Hamilmonian flow in any dimension \cite{BCI23}. 

Yet, these examples are very fragile, while there are abundant examples among dissipative systems \cite{Ja81,BC2}  which display behavior similar to what is numericaly observed in the standard map. Yoccoz' strong regularity program \cite{Yolecture1,BY19} proposed to  set up such proofs in a formalism which enables to deduce hyperbolic dynamical properties from topological and combinatorial properties, in order to show the abundance of  attactors of the higher and higher Hausforff dimension. See  \cite{Y19,PY09,berhen} for the first steps of this program.\\

\paragraph{\bf Growth of the number of periodic points}  For every  polynomial map, the number of isolated $n$-periodic points grows at most exponentially fast with $n$. This is also the case of a $C^\infty$-dense set $\mathcal D$ in $\Diff^\infty(M)$, see \cite{AM65,K99}. This led Smale \cite{Sm} and Arnold \cite{Ar00} to wonder if a topologically generic dynamical systems has its number of periodic points growing at most exponentially. This turns out to be true in dimension 1 by \cite{MdMvS92} and in some real-analytic settings \cite{KH07,KH072}.
  
On the other hand, whenever $2\le r<\infty$ and $\dim M\ge 2$, Kaloshin proved in \cite{K00} that a topologically generic diffeomorphism in an open subset of $\Diff^r(M)$ has the number of periodic points growing super-exponentially. The conservative counterparts of this result were obtained in \cite{KS06,GTS07}. Kaloshin's approach employed periodic bands (lines of periodic orbits) emerging out of a homoclinic tangency. For classes of systems without homoclinic tangencies, the super-exponential growth was established in \cite{BerPer,AST16,AST20}. 

For the typicality in a stronger sense, the following conjecture was proposed \footnote{Many Arnold's problems are related to this question \cite[Problems 1992-13, 1994-48, 1992-14]{Ar00}.}:
\begin{conj}[Pb. 1994-47  \cite{Ar00}]\label{problem2} 
The number of periodic points of a mapping of class $C^\infty$  grows almost always not faster than some exponential function of the period. Here ``almost always" means ``for almost all (in the sense of the Lebesgue
measure) parameter values in each typical family of mappings depending on
sufficiently many parameters.''
\end{conj}

In  \cite{As16}, Asaoka showed that among surface symplectomorphisms with KAM circles, super-exponential growth of the number of periodic points is Kolmogorov typical (also for real-analytic maps). The Kolmogorov typicality of the super-exponential growth for a class of (non-conservative) $C^r$-generic diffeomorphisms with normally-hyperbolic invariant circles
was proved in \cite{BerPer}. Note that our Theorem \ref{Main explicit per} generalizes results of \cite{As16} to higher dimensions and invalidates the symplectic counterpart of Conjecture \ref{problem2} in this context. 
   
\section{Elliptic points and KAM tori} \label{accum KAM}
Consider a $C^\infty$ symplectomorphism $f$ with a $q$-periodic point $O$. Assume $O$ is \bfemph{elliptic},
i.e., the eigenvalues 
$\lambda_1,\dots,\lambda_{2n}$ of the derivative $D f^q$ at $O$ are simple and lie on the unite circle:
$$\lambda_k = e^{i\hat \alpha_k}, \quad \lambda_{k+n}=e^{-i\hat\alpha_k},\qquad 1\leq k \leq n,
$$
for a set of $n$ different reals $\hat\alpha_k\in (0,\pi)$. The next assumption is that $O$ has \bfemph{no strong resonances}, meaning that
\[
\left\{\begin{array}{c} \lambda_{j_1} \lambda_{j_2}\lambda_{j_3} \neq 1 \\
\lambda_{j_1} \lambda_{j_2}\lambda_{j_3}\lambda_{j_4} = 1\Rightarrow 
\text{
two of these eigenvalues are the conjugates of the two others.
}\end{array}\right. \]

\begin{theorem*}[Birkhoff \cite{Bir27}] 
In a neighborhood of the elliptic point $O$ that has no strong resonances, one can introduce symplectic local coordinates $z=(z_k)_{1\le k\le n}\in \C^{n}\approx \R^{2n}$ such that 
\begin{equation}\label{bnfel}
f^q(z)= \left(z_k\cdot \exp  [ i \hat\alpha_k +i \sum_m \omega_{k, m}\cdot  |z_m|^2 ]\right)_{1\le k\le n}  +\xi(z,z^*),
\end{equation}
for some set of real numbers $\omega_{k,m}$, $1\leq k \leq n$,  $1\leq m \leq n$, and a $C^\infty$ function $\xi$ vanishing at zero along with all derivatives up to order 3. The coordinate transformation depends smoothly on the map $f$ in $C^\infty$. 
\end{theorem*}
 
\begin{definition}\label{def_non_deg} An elliptic periodic point is \bfemph{non-degenerate} if there are no  strong resonances and the \bfemph{twist condition} holds: the matrix $ [\omega_{k, m}]_{1\le k, m\le n}$ is invertible. 
\end{definition}

The non-degeneracy condition is open: 
at small symplectic perturbations of $f$,  the elliptic points persist  and
both $\hat\alpha_k$ and $\omega_{k,m}$ depend on $f$ continuously. Also, given an elliptic point, it can be made non-degenerate by an arbitrarily small symplectic perturbation of $f$. 

In the proof of Theorem \ref{main1} we use the following fact, which is one of the basic applications of the KAM-theory: 
\begin{proposition}\label{accumulation}  Any non-degenerate elliptic point of a smooth symplectomorphism is accumulated by non-degenerate KAM tori.
\end{proposition}

We provide a proof below, for the sake of completeness (the analytic case has been shown in 
\cite{El88}). First, let us remind the definitions.

A vector $\alpha\in \R^n$ is \bfemph{Diophantine} if there exist 
$\gamma>0$ and $\tau>0$ such that, for every $k\in \Z^{n}$ and $l\in  \Z$,
\begin{equation}\label{dio}
 |l +<k, \alpha>|>  \gamma \cdot \|k\|^{-n-\tau}.
\end{equation}
The corresponding map $R_\alpha: \theta\in \T^n\mapsto \theta+\alpha\in \T^n$ is called a  \bfemph{Diophantine rotation}.
  
A Lagrangian $n$-torus $T\subset M$ is a \bfemph{KAM torus} if it is invariant with respect to an iteration of $f$, which acts on $T$ as a Diophantine rotation. Namely, there exist an integer $q\ge 1$, a $C^\infty$ embedding $h$ of $\T^n=\R^n/\Z^n$ into $M$, and 
a Diophantine $\alpha\in \R^n$ such that $T=h(\T^n)$ and
\[f^q\circ h= h\circ R_\alpha.\]
\begin{remark}\label{KAM Tori disjoint}
It two KAM tori intersect, they must be equal. Indeed their intersection is closed, and is also dense, as the dynamics on each torus are minimal. 
\end{remark}
According to \cite{Bir36} (see \cref{Birkhoff} below), near a KAM-torus $T$ one can introduce symplectic $C^\infty$ coordinates $(\theta,r)\in \T^n\times \R^n$ such that  
\[f^q(\theta, r) =(\theta + \alpha + \Omega r, r)+O(r^2),\] 
where $\Omega$ is a symmetric $(n\times n)$ matrix.
\begin{definition} The KAM torus is called \bfemph{non-degenerate} if $\det \Omega \neq 0$.\end{definition}

The KAM theory ensures the \bfemph{persistence} of non-degenerate KAM tori: 
\begin{theorem}[\cite{Mo62,Bo84}] \label{KAM classique}
Let $T=h(\T^n)$ be a non-degenerate KAM torus with a Diophantine frequency vector $\alpha$ for some exact $f\in \Symp^\infty(M)$. Then, any exact symplectomorphism
$\hat f$ which is $C^\infty$-close to $f$ has a uniquely defined non-degenerate KAM torus close to $T$, with the same frequency vector $\alpha$. 
In particular, if $(f_a)_a$ is a $C^\infty$ family of symplectomorphisms such that $f_0=f$, then for all small $a$ there exists a $C^\infty$ embedding 
$h_a: \T^n \hookrightarrow M$ such that $h_0=h$ and, for all small $a$, the torus $h_a(\T^n)$ is the uniquely defined non-degenerate KAM torus for $f_a$ satisfying 
$f_a^q\circ h_a = h_a\circ R_\alpha$. Moreover,  the family $(h_a)_a$ is $C^\infty$.
\end{theorem}

\begin{proof}[Proof of \cref{accumulation}]
Choose symplectic coordinates such that $f^q$ near the elliptic point $O$ has the 
form (\ref{bnfel}). Choose any small $\varepsilon>0$, and make the following (conformally symplectic) coordinate change in a small neighborhood of the $n$-torus $\{|z_k|=\varepsilon\}_{k=1,\ldots,n}$:
\[ z_k= \varepsilon \sqrt{1 +\varepsilon^{1/3} r_k} \cdot \exp[i \theta_k],\]
where $(\theta,r)=(\theta_1, \dots , \theta_n ,r_1,\dots,  r_n) \in  \T^n \times \R^n$.
The substitution into (\ref{bnfel}) gives
\[f^q(z)=
 \left(\varepsilon \sqrt{1 +\varepsilon^{1/3} r_k} \cdot  \exp[i (\theta_k+\hat\alpha_k+ \sum_m \omega_{k, m}\cdot  (\epsilon^2 +\epsilon^{7/3} r_m) )]\right)_{1\le k\le n}  
 + \epsilon ^4 \cdot   \psi_\epsilon(\theta,r), \]
for a $C^\infty$ function $\psi_\epsilon: \T^n\times \B^n \to  \C^n$ which 
is uniformly $C^\infty$-bounded as $\eps\to 0$ (i.e., it
stays uniformly $C^\rho$-bounded as $\epsilon\to 0$, for any given 
$\rho\geq 1$).
Hence, in the coordinates $(\theta,r)$, the map $f^q$ has the form
 \[(\theta, r) \mapsto \left( \theta_k + \tilde\alpha_k +  \varepsilon^{7/3}\sum_m \omega_{k, m}\cdot r_m,\; 
r_k\right)_{1\le k\le n}   +  \varepsilon^{8/3}\cdot  \phi_{\epsilon}(\theta,r), \]
where $\tilde\alpha_k =\hat\alpha_k +  \varepsilon^{2}\sum_m \omega_{k, m} $, and $\phi_\epsilon$ is uniformly $C^\infty$ bounded as $\eps\to 0$, on any compact set of $(\theta,r)$. 

We now assume that $\epsilon= s^{-3/7}$ for a large integer $s$. This gives:
\begin{equation}\label{fqcl}
f^q: (\theta, r) \mapsto \left( \theta_k + \tilde\alpha_k + s^{-1}\sum_m \omega_{k, m}\cdot r_m,\; 
r_k  \right)_{1\le k\le n} + s^{-8/7}\cdot \tilde \phi_s(\theta,r),
\end{equation}
where  the functions $\phi_s$ stay uniformly $C^\infty$-bounded as $s\to+\infty$, on any compact set of $(\theta,r)$.

Iterating $s$ times, we obtain from (\ref{fqcl})
\[f^{sq}: (\theta, r) \mapsto (\theta + s(\tilde\alpha_k)_{1\le k\le n}  +    [\omega_{k, m}] \cdot r,\;  r)  +  s^{-1/7}\phi_q(\theta,r),\]
where the functions $\phi_s$ are uniformly $C^\infty$ bounded as $s\to+\infty$, on any compact set of $(\theta,r)$.
Thus, taking a sequence of the integers $s$ going to infinity such that the vector $(s\tilde\alpha_k \ \mod 1)_{1\le k \le n}$ tends to a vector  
$\alpha \in \R^n$, we see that the maps $f^{sq}$ get $C^\infty$-close to the integrable map
\begin{equation} \label{themap} (\theta, r ) \mapsto (\theta + \alpha  +    [\omega_{k, m}] \cdot  r,\;  r),\end{equation} 
which satisfies the twist condition $\det [\omega_{k, m}] \neq 0$. Let $\alpha_o\in \R^n$  be a Diophantine vector  sufficiently close to $\alpha$ such that there exists $r_o\in \B^n$ satisfying:
\[ \alpha  +    [\omega_{k, m}] \cdot  r_o= \alpha_o .\]
  By definition,  the map  (\ref{themap})   has an invariant non-degenerate KAM-torus $\{r=r_o\}$ of frequency $\alpha_o$. By the persistence theorem
(\cref{KAM classique}), each of the maps $f^{sq}$ with sufficiently large $s$ and such that $(s\tilde\alpha_k \ \mod 1)_{1\le k \le n}$ is close to  $\alpha$ have an invariant non-degenerate KAM-torus $T_s$ close to $r=r_o$, with the frequency vector $\alpha_o$.

Note that $T_s$ is invariant with respect to $f^q$. Indeed, by (\ref{fqcl}) the torus $f^q s$ is close to $r=r_0$; this torus is invariant with respect to $f^{sq}$ and has the same frequency vector $\alpha_o$ as $T_s$, so $f^qT_s=T_s$ by uniqueness.
\end{proof}

It is useful to represent the symplectic map near a KAM torus with the help of a generating function. Let $(\theta,r)\in \T^n \times \B^n$ be
symplectic coordinates near a KAM torus $T: \{r=0\}$. We assume that $T$ is $q$-periodic, i.e., $f^qT=T$ for some $q\geq 1$ and $f^jT\cap T =\emptyset$ for $1\leq j <q$. Denote 
$f^q(\theta,r)=(\bar \theta, \bar r)$.
As $f^q$ is exact symplectic, $\bar r d\bar \theta- r d \theta$ is exact and so is $-\bar r d\bar \theta+ r d \theta + d(r (\bar \theta-  \theta) )$. Hence there exists a function $S$ such that:
\[ dS= -\bar r d\bar \theta+ r d \theta + d(r (\bar \theta-  \theta) )
= (r-\bar r) d\bar \theta +(\bar \theta-  \theta) dr
\]
 As $f^q(\theta, 0) = (\theta + \alpha, 0)$ for a constant 
$\alpha$, the derivative $\partial_\theta\bar\theta$ is invertible at small
$r$. Hence, by the implicit function theorem, one can consider
$(\theta, \bar r)$ as a function of $(\bar\theta, r)$. 
The function $S(\bar\theta,r)$ is called a \bfemph{generating function}. 
By the construction,
$f(\theta, r) = (\bar \theta, \bar r)$ if and only if
\begin{equation}\label{sd0}
\theta = \bar\theta-\partial_{r} S (\bar\theta, r)  \qand  \bar r = r-\partial_{\bar \theta}  S (\bar\theta, r).
\end{equation}

 Since $\bar\theta= \theta+\alpha +O(r)$, we have here $S=\alpha r + O(r^2)$. Note that when $f=f_a$ depends smoothly on a parameter $a\in \sA$,  we can choose $S=S_a$ depending smoothly on $a$. 

A theorem by Birkhoff establishes \cite{Bir36} that dynamics near a KAM-torus  are integrable up to flat terms. We will use the following parametric version of this result:
\begin{proposition}\label{Birkhoff} Let $(f_a)_a$ be a $C^\infty$ family of symplectomorphisms and let $(T_a)_a$ be a $C^\infty$ family of $q$-periodic KAM tori for the maps $f_a$ with the Diophantine vector 
$\alpha$ independent of $a$. Then, for every integer 
$D\geq 0$, in a small neighborhood of the torus $T_a$ there exist $C^\infty$-smooth symplectic coordinates $(\theta,r)\in \T^n \times \R^n$ (depending $C^\infty$ smoothly on $a$) such that $f^q_a(\theta,r)=(\bar\theta,\bar r)$ 
if and only if
\begin{equation}\label{nfcrb}
\bar\theta = \theta + \alpha + \nabla Q_a(r) +  \partial_r s_a(\bar \theta, r),\qquad
\bar r = r - \partial_{\bar \theta} s_a(\bar \theta, r),
\end{equation}
where $(Q_a)_a$ is a $C^\infty$ family of degree-$(D+1)$ polynomials satisfying $\nabla Q_a(0)=0$, and  $(s_a)_a$ is a $C^\infty$ family of functions that vanish at $r=0$ along with all derivatives with respect to $r$ up to the order $(D+1)$.
\end{proposition}
\begin{proof}
We proceed by induction. Our goal is to show that for any $D\geq 0$ a symplectic coordinate transformation can keep the map $f^q_a$ in the form (\ref{sd0}) while bringing the generating function $S_a$ to the form
\begin{equation}\label{sd1}
S_a(\bar\theta, r) =  \alpha r + Q_{a,D}(r)  + s_{a,D}(\bar\theta,r),
\end{equation}
where $Q_{a,D}$ is a polynomial of $r$ of degree $(D+1)$ (with the coefficients $C^\infty$-smoothly dependent on $a$), and $s_{a,D}$ has all its derivatives with respect to $r$ vanish at $r=0$ up to the order $(D+1)$.

This is true at $D=0$ with $Q_{a,0}=0$. So, assume that the generating function is in the form (\ref{sd1}) at some $D\geq 0$ and let us show that it can be brought to the same form at $(D+1)$, by a symplectic transformation.

Rewrite (\ref{sd1}) as 
\begin{equation}\label{sd3}
S(\bar\theta, r) =  \alpha r + Q_{a,D}(r) + \Omega_{a}(r) + \hat\Omega_{a}(\bar \theta,r) + \hat s_a(\bar \theta, r),
\end{equation}
where $\Omega_{a}(r)$ is a homogeneous polynomial of $r$ of degree $(D+2)$ and $\hat\Omega(\bar \theta,r)$ is a degree-$(D+2)$ homogeneous polynomial of $r$ with coefficients which, as functions of $\bar\theta$, have zero mean;
the function $\hat s_a(\bar \theta, r)$ has all derivatives with respect to $r$ vanish at $r=0$ up to the order $(D+2)$. 
As $\bar\theta = \theta + \alpha + O(r)$ and $\hat\Omega_{a}$ is of degree $D+2$ in $r$, it holds:
\[  \hat\Omega_{a}(\bar \theta,r) = \hat\Omega_{a}(  \theta+\alpha,r)+ 
o(r^{D+2}).\] 
Then   by  (\ref{sd0}),(\ref{sd3}) we have that:
\begin{equation}\label{sd2}
\bar\theta = \theta + \alpha + O(r), \qquad
\bar r = r - \partial_{ \theta} \hat \Omega_a( \theta+\alpha,r) + 
o(r^{D+2}).
\end{equation}

The function $\hat\Omega_a$ has zero mean, so we can write its Fourier expansion as

$$\hat\Omega_a(\theta +\alpha,r)= \sum_{m\in\Z^n, m\neq 0} \hat\Omega_m(r;a) e^{i m\theta},$$
where $\hat\Omega_m(r;a)$ are homogeneous polynomials of $r$ of degree $(D+2)$ with $a$-dependent ($C^\infty$) coefficients, which decay to zero faster than any power of $\|m\|$ as $m\to \infty$. 

Define a $C^\infty$ function $\hat Q_a(\theta,r;a)$ (a homogeneous polynomial of degree $(D+2)$ in $r$) by the rule
$$\hat Q_a(\theta,r;a)=\sum_{m\in\Z^n, m\neq 0} \frac{e^{im\alpha}}{1-e^{im\alpha}}\hat\Omega_m(r;a) e^{i m \theta};$$
the series is convergent since $\alpha$ is Diophantine.
This function satisfies the identity
\[\hat Q_a(\theta+\alpha, r) 
=\hat Q_a(\theta, r) - \hat\Omega_a(\theta+\alpha, r). \]
One sees from this formula and (\ref{sd2}) that
$$ \bar r - \partial_{\bar\theta} \hat Q_a( \theta+\alpha,r)
=
\bar r + \partial_{\theta}\hat\Omega_a(\theta+\alpha, r) - \partial_{\theta} \hat Q_a(\theta, r)=
r - \partial_{\theta} \hat Q_a(\theta, r) + o(r^{D+2}).$$
Using again \eqref{sd2} and the  homogenity of the polynomials, we have
$  \partial_{\bar\theta} \hat Q_a(\bar\theta,\bar r) =  \partial_{\bar\theta} \hat Q_a( \theta+\alpha,r) +  o(r^{D+2})$ and so
$$ \bar r - \partial_{\bar\theta} \hat Q_a( \bar \theta ,\bar r)
=
r - \partial_{\theta} \hat Q_a(\theta, r) + o(r^{D+2}).$$

This means that after the symplectic coordinate transformation
$$(\theta,r)\mapsto \Phi (\theta,r;a) = (\theta + \partial_r \hat Q_a(\theta,r) + o(r^{D+1}), r - \partial_\theta \hat Q_a(\theta, r) + o(r^{D+2})),$$
where $\Phi$ denotes here the time-1 map of the flow defined by the Hamilton function $\hat Q_a$,
we get
$$\bar r - r = o(r^{D+2}).$$
Thus, $\bar r - r$ is flat up to degree $(D+2)$ in the new coordinates. Therefore, the corresponding generating function cannot have $\theta$-dependent terms of degree $(D+2)$ or lower, i.e.
$$S(\bar\theta, r) = \alpha r + Q_{a,D+1}(r)  + o(r^{D+2}),$$
as required (one may check that the degree-$(D+2)$ polynomial $Q_{a,D+1}(r)$ equals to
$Q_{a,D}(r) + \Omega_a(r)$, see (\ref{sd3})). Importantly note that the latter coordinate change depends smoothly on $a$.
\end{proof}
We will call the map (\ref{nfcrb}) the Birkhoff normal form, and the corresponding coordinates $(\theta,r)$ the normalizing coordinates.

\section{Proof of Theorem \ref{main1}}\label{proof main1}
 Let $\sA$ be a  compact manifold and let  $(f_a)_{a \in \sA}\in \E(M)_\sA$. We shall prove that for every $D\geq 1$ there exists a $C^D$-close family $(\tilde f_a)_{a \in \sA}$ for which there is an open covering $(\sA_i)_{1\le i\le I}$  of $\sA$ and, for every $1\le i\le I$, there exist  a nonempty open subset $U_i$ and an integer $N_i\ge 1$ such that  $f^{N_i}_a|_{U_i}=id$ for every $a\in \sA_i$.
\medskip
 
The following shown below will be used several times:
\begin{lemma} \label{lemprt}
Let $V\subset M $ be an open subset of a symplectic manifold 
$(M,\omega)$.  Let $V'\Subset V$ be a relatively compact subset.  
For every $m\ge 0$, let  $g_{\sA,m}:=(g_{a,m})_{a\in \sA} $ be a family of exact symplectic maps $g_{a,m}:  V\to M$ such that the sequence 
$(g_{\sA,m})_{m} $ converges to the family constantly equal to the canonical inclusion  $id_V: V\!\hookrightarrow\!M$ (in the $C^\infty$ compact open topology). 
Then there exists a convergent to $(id_V)_{a\in \sA} $ sequence of families  $\tilde g_{\sA,m}:=(\tilde g_{a,m})_{a\in \sA} $ of exact symplectic maps which coincide with $g_{a,m}$ on $V'$ and are compactly supported. 
\end{lemma}

\begin{proof} 
%
The diagonal $\Delta$ of $M\times M$ is a Lagrangian submanifold of the symplectic manifold $(M,\omega)\times (M,-\omega)$, therefore, by the Weinstein tubular neighborhood theorem, a neighborhood of $\Delta$ is $C^\infty$ symplectomorphic to a neighborhood of $\Delta\times \{0\}$ in the cotangent bundle $T^*\Delta$ of $\Delta$. The graph of the canonical inclusion $id_V$ is sent by this symplectomorphism to  $V\times \{0\}$, so the graph of each $g_{a,m}$ (with $m$ large) is identified with a section $\sigma_{a,m}$ of $T^*\Delta$ estricted to an open subset of $\Delta$. Observe that the family $(\sigma_{a,m})_{a\in \sA}$ is smooth and converges to the zero section
as $m\to\infty$. Also by \cite[Lem. 10.2.8]{MDS17}, as each $g_{a,m}$ is exact, the section $\sigma_{a,m}$ is exact: there exists $\phi_{a,m} \in C^\infty(\Delta, \R)$  such that  $\sigma_{a,m}= d \phi_{a,m}$, 
and we can chose the families $(\phi_{a,m})_{\sa\in \sA}$ smooth and converging to $0$ when $m\to\infty$. 
 
 Now take a function $\rho\in C^\infty(V, \R)$,  compactly supported in $V$, which is equal to $1$ in a neighborhood of $cl(V')$. Observe that when $m$ is large, for every $a$,  the domain of $\phi_{a,m}$ contains $\Delta \cap (\mathrm{Supp}\, \rho)^2$. 
   Hence for $m$ large enough, we can define  $\tilde \phi_{a,m} \in C^\infty(\Delta,\R)$   as equal to  the function  with support in $\Delta \cap (\mathrm{Supp}\, \rho)^2$ and such that $ \tilde \phi_{a,m}(y)= \rho(x)\cdot   \phi_{a,m}(y)$ for every $y=(x,x)\in \Delta \cap   (\mathrm{Supp}\, \rho)^2$. 
Then by  \cite[Lem. 10.2.8]{MDS17}, the exact section ~$d \tilde \phi_{a,m}$ is the graph of an exact symplectomorphism $\tilde g_{a,m}$. Note that it satisfies the requested properties.
\end{proof} 
 
For every $a\in \sA$, by \cref{accumulation}, the map $f_a$ has a non-degerate KAM torus. Then, by \cref{KAM classique}, this KAM torus persists with the same frequency vector, for all sufficiently close $a$. 
Hence, by \cref{Birkhoff}, there is an open  covering $(\sA_i)_i$ of $\sA$ such that for every $i$, for some integer $q_i\geq 1$ there exists a $q_i$-periodic
KAM-torus $T_{a,i}$ with the Diophantine frequency vector $\alpha_i$, such that, for some choice ($C^\infty$ dependent on $a\in \sA_i'$)  of symplectic coordinates $(\theta,r) \in \T^n\times \R^n$ in a neighborhood $V_{a,i}$ of $T_{a,i}:\{r=0\}$,
the map $f_a^{q_i}|_{V_{a,i}}:(\theta,r)\mapsto (\bar\theta,\bar r)$ is represented in the form \eqref{nfcrb}, with $\alpha=\alpha_i$ and $i$-dependent functions $Q_a$ and $s_a$. 
So, the map $f_a^{q_i}$ in $V_{a,i}\cap f_a^{-q_i}(V_{a,i})$ is defined by the generating function 
$$S_{a, i}= \alpha_i r + Q_{a,i}(r) + s_{a,i}(\bar\theta,r).$$
Since $T_{a,i}$ is non-degenerate, we have 
\begin{equation}\label{nfgR}
\det (\nabla^2 Q_{a,i}(0))\neq 0.
\end{equation}

By compactness of $\sA$, we can assume that the covering $(\sA_i)_i$ is finite, $1\leq i \leq I$. 
Note that, for every $a$, the tori $T_{a,i_1}$ and $T_{a,i_2}$ can be assumed disjoint for $i_1\neq i_2$.
Indeed, otherwise  $T_{a,i_1}=T_{a,i_2}$ for some $a$
and $i_1\neq i_2$ by \cref{KAM Tori disjoint}, and we can modify the covering $(\sA_i)_i$ by replacing the elements $\sA_{i_1}$ and $\sA_{i_2}$ by their union (the normalizing coordinates near $T_{a,i}$ which were chosen separately for $a\in \sA_{i_1}$ and $a\in \sA_{i_2}$ need to be replaced then by the normalizing coordinates that are defined, and are $C^\infty$ smooth, for all $a\in \sA_{i_1}\cup \sA_{i_2}$).  

By shrinking $V_{a,i}$ and invoking \cref{KAM Tori disjoint}, we can always assume that the consecutive images $f_a^{j}V_{a,i}$, $0\leq j < q_{i}$, $1\leq i\leq I$, are disjoint for all $a$. Let $\kappa(u)$ be a $C^\infty$ function that vanishes for all $|u|>2$ and equals to $1$ for all $|u|<1$. Note that 
$$S_{\delta, a, i} = \alpha_i r + Q_{a,i}(r) + (1-\kappa(\delta^{-1}\|r\|))\;s_{a,i}(\bar\theta,r),$$
form a $C^{D}$-small family perturbation of the family of generating functions $(S_{a,i})_{a\in \sA_i}$
if $\delta$ is small enough (we use here that $s_{a,i}$ vanishes at $r=0$ along with all derivatives with respect to $r$ up to the order $(D+1)$). 
This defines a family of exact symplectomorphisms  which is $C^D$-close to the identity. 
Hence, by \cref{lemprt}, there is a $C^D$-small perturbation $(\tilde f_a)_{a \in \sA}$ of $(f_a)_{a\in \sA}$  such that for every $1\leq i \leq I$, for slightly shrunk sets $\sA_i$ and neighborhoods $V_{a,i}$, the map $\tilde f_a^{q_i}$ in $V_{a,i}\cap f_a^{-q_i}(V_{a,i})$ is defined by the generating function $S_{\delta,i}$. Shrinking $V_{a,i}$ to $V_{a,i}\cap \{\|r\|<\delta\}$, we obtain that the maps 
$\tilde f_a^{q_i}$
for the perturbed family $\tilde f_a$ are defined by the generating functions $\tilde S_i=\alpha_i r + Q_{a,i}(r)$, i.e.,
these maps become integrable:
\[\tilde f_a^{q_i}|_{V_{a,i}}(\theta,r) = (\theta + \alpha_i  + \nabla Q_{a,i}(r), r).\]

Next, by appealing to \cref{lemprt} again, we can further perturb the family $(\tilde f_a)_{a\in \sA}$ such that the
frequency vectors become rational, i.e., we may assume $\alpha_i\in {\mathbb Q}^n$ for all $1\leq i \leq I$.
Let $p_i$ be the minimal positive integer such that $p_i\cdot \alpha_i\in \N^n$. Then 
\begin{equation}\label{integrma}
\tilde f_a^{p_iq_i}|_{V_{a,i}}(\theta,r) = (\theta  + p_i\nabla Q_{a,i}(r), r).
\end{equation}

By construction, $O_{a,i}=(0,0)$ is a point of the least period $p_iq_i$ for the map $\tilde f_a$.
As $r=0$ is a non-degenerate critical point of $Q_{a,i}$ (see (\ref{nfgR})), by Morse lemma there exists
a diffeomorphism $r\mapsto Y_{a,i}(r)$ of a neighborhood of zero in $\R^n$ such that
$$p_i Q_{a,i}(r) = \frac{1}{2} <Y_{a,i}(r), R_i Y_{a,i}(r)>$$
where $R_i$ is a real diagonal $(n\times n)$-matrix with non-zero diagonal entries (equal to $\pm 1$).
By \cite[P. 502]{Ho85}, the family  $(Y_{a,i})_{a\in \sA_i}$ can be chosen $C^\infty$ smooth.

Let us introduce symplectic coordinates $(x,y)\in B^n\times B^n$ in a small neighborhood of the periodic point $O_{a,i}$
by the rule
$$\theta = (\partial_r Y_{a,i}(r))^\top x , \qquad y =Y_{a,i}(r),$$
where $^\top$ denotes the transpose of a matrix.
To see that this coordinate  transformation is symplectic, note that it is defined by a generating function: $y=\partial_x S, \; \theta = \partial_r S$, with $S=<x,Y_{a,i}(r)>$. 
 
In these new local coordinates, the return map (\ref{integrma}) near the periodic point $O_{a,i}$ takes the following form:
$$\tilde f_a^{p_i q_i}: (x, y)\mapsto (x+ R_i y, y),$$
for all $a\in \sA_i$. Since $p_iq_i$ is the least period of $O_{a,i}$, \cref{lemprt} gives that, by 
a $C^\infty$-small perturbation of the family $(\tilde f_a)_{a\in \sA}$, the return maps near the periodic points $O_{a,i}$
can be made equal to
$$\tilde f_a^{p_i q_i}: (x, y)\mapsto ((1-\eps) x+ R_i y, -\eps R_i^{-1}x + y),$$
for any small $\eps>0$, the same for all $a\in \sA$ and $1\leq i\leq I$; recall that $R_i$ are diagonal matrices.

This is a symplectic linear map with complex eigenvalues equal to $1$ in the absolute value. Moreover, for a neat choice of 
$\epsilon$, they are roots of $1$, the same for all $a\in \sA'_i$: take $\eps = 2(1-\cos \frac{2\pi}{s})$. So, 
for the perturbed family $\tilde f_a$ and $N_i=sp_iq_i$, we have that $\tilde f_a^{N_i}=id$ in
some neighborhood $U_{a,i}$ of $O_{a,i}$.
\qed

\section{Emergence}
\subsection{Definition of the emergence order}\label{section def emergence}
Given a volume-preserving map $f$ of a bounded manifold $M$ endowed with a normalized volume form $\leb$, the statistical behavior of a point $z$ is given by the empirical measure:
\[\mathsf e(z):= \lim_{n\to \infty}  \frac1n \sum_{k=0}^{n-1} \delta_{f^k(z)}\]
(the limit exists for $\leb$-a.e. $z$ by the Birkhoff-Khinchin theorem). This is a probability measure on $M$. The map $z \mapsto \mathsf e(z)$ is a a measurable map from $M$ into the space $\cM_1(M)$ of probability  measures on $M$.  The ergodic decomposition of $f$ is the push-forward $\mathsf e_*\leb$. It is a probability measure on the space
$\cM_1(M)$. Basically, it describes the distribution of statistical behaviors of the orbits of $f$. For instance, if $f$ displays a finite number of physical measures $(\mu_i)_{1\le i \le I}$ whose basins $(B_i)_i$  cover $\leb$-a.e. points of $M$,  then
\[\mathsf e_*\leb = \sum_{i=1}^I\leb(B_i)\cdot  \delta_{\mu_i}\]  
The emergence quantifies how far we are from this situation. To define it, we endow the space  $\cM_1(M)$ with the Kantorovich-Wasserstein distance $d_{KW}$; this induces the weak-$*$ topology. To recall the definition, 
given probability measures $\mu$ and $\nu$, we regard transport plans between  $\mu$ to $\nu$: these are probability measures on $M^2$ which push forward to $\nu$ via the first coordinate projection and to $\mu $ via the second coordinate projection. Then, with $\Pi(\nu, \mu)$ being the space of the transport plans between these measures, we define:
\[d_{KW}(\nu, \mu)=  \inf_{\pi\in \Pi(\nu, \mu)} \int_{(x,y)\in  M^2}   d(x,y) \,  d\pi, \]
where $d(x,y)$ is the distance in $M$.
\begin{definition}\label{def emergenece}
The \bfemph{Emergence of $f$ at scale $\epsilon$} is the minimum number $N$ of probability measures 
$(\mu_i)_{1\le i\le N}$ such that:
\[\int_M \min_{1\le i\le N}  d_{KW} (\mathsf e(z),\mu_i) d\leb < \epsilon.\]
We denote it by $\sE(\epsilon):=N$. 
We say that the dynamics $f$ has \bfemph{emergence of order $\alpha\ge 0$} if 
\[\limsup_{\epsilon \to 0} - \frac{\log \log \sE(\epsilon) }{\log \epsilon} = \alpha.\]
\end{definition}
It was shown in \cite{BB21} that the order of emergence is at most $\dim M$, for the covering number $\cN(\epsilon)$
of $\cM_1(M)$ has order $\dim M$:
\[\lim_{\epsilon \to 0} - \frac{\log \log \cN(\epsilon) }{\log \epsilon} = \dim M.\]
Hence, we say that the order of the emergence is \emph{maximal} if it equals $\dim M$.  

Given an invariant open subset $E\subset M$ of positive Lebesgue measure, we can consider the normalized measure $\leb_E$ of $\leb$ restricted to $E$ and consider likewise the emergence of $(f|_{E} , \leb_E)$.  An immediate consequence of \cite[Lem. 3.18]{BB21} is:
\begin{lemma}\label{relation emer rest}
The emergence of $f$ at scale $\leb(E)\cdot \epsilon$ is at least the emergence of $f|_E$ at scale $\epsilon$. 
\end{lemma}

\subsection{Emergence of the maximal order is a localizable property}\label{emergence proof}
The aim of this section is to prove \cref{example of localizable prop3} for symplectomorphisms of a compact manifold $(M, \omega)$ of any dimension $2n\ge 2$. 
We notice that  \cref{example of localizable prop3} is a consequence of  the following:
\begin{lemma}\label{CS1}
Take any $\tau_0<2n$ and any $f\in \Symp^\infty(M)$ which displays a domain $U\subset M$ filled by periodic points of the same period. Then, for every $\epsilon_0>0$ there exist
$\epsilon<\epsilon_0$ and an open set $\cU$ near $f$ in $\Symp^\infty(M)$ such that for every $\tilde f\in \cU$, the emergence at scale $\epsilon$ of $f$ is at least $exp(\epsilon^{-\tau_0})$. 
\end{lemma}

Let $p\ge 1$ be the minimal integer such that $f^p|_{U}= id$.  By shrinking $U$, if necessary, we can assume that $f^k(U)$ intersects $U$ only if $k$ is a multiple of $p$. 

\begin{lemma}\label{CS2}
There exists $\lambda>0$ such that for any $\tau<2n$ and $\epsilon_0>0$, there exist $\epsilon<\epsilon_0$ and a nonempty open set $\cal N_{\tau,\epsilon}\subset \Symp^\infty(U\! \hookrightarrow\!  M)$ near the canonical inclusion $U\hookrightarrow M$  such that every $g\in \cal N_{\tau,\epsilon}$  leaves invariant a subset $U'\subset U$ of measure $\ge \lambda$, such that the emergence at scale $\epsilon$ of $g|_{U'}$ is at least $\exp(\epsilon^{-\tau})$ (for $U'$ endowed with the normalized measure $\leb|_{U'}$).
\end{lemma}
We postpone the proof of this lemma for a short while, and show, first, that \cref{CS2} implies  \cref{CS1}, hence \cref{example of localizable prop3}. Indeed, by \cref{lemprt}, for  
any compactly supported perturbation $g$ of the identity in $\Symp_c^\infty (U)$,  there exists a perturbation $\tilde f$  of $f$, localized in $U$,  such that $\tilde f^p|_U= g$. Thus, for any  $\tau < 2n$ and suitable $\epsilon>0$, the set $\hat \cN_{\tau,\epsilon} \subset \Symp^\infty (M)$ of maps $\tilde f$ such that $\tilde f^p|_U$ belongs to $\cN_{\tau,\epsilon}$ is a nonempty open set that lies near $f$.  So, it remains to show that for every $\tau_0$ and $\epsilon_0>0$, there exists  $\tau<2n$ and $\epsilon<\epsilon_0$ such that $\hat \cN_{\tau,\epsilon}$ is formed by maps whose emergence at scale $\epsilon$ is greater than $\exp(\epsilon^{-\tau_0})$. 
 
Let $\epsilon_0$ and $\tau_0$ be as in the statement of \cref{CS1}.  Let $\tau:= n+\tfrac12 \tau_0$; note that $2n>\tau>\tau_0$. 
Let $\epsilon<\epsilon_0$ satisfy \cref{CS2} and be sufficiently small, so
\[\lambda^{-\tau} \cdot \epsilon^{\tau_0-\tau}>1.\]
Let $\tilde f\in \tilde \cN_{\tau,\epsilon}$.  As $\leb|_U \ge \lambda\cdot  \leb$, it follows by \cref{relation emer rest} that the emergence of $\tilde f$ at scale $\epsilon$ is at least
  \[\exp((\lambda \epsilon)^{-\tau})\ge  
  \exp( \lambda^{-\tau}\cdot  \epsilon^{-\tau +\tau_0}
  \cdot \epsilon ^{-\tau_0}
  )\ge  
  \exp( \epsilon^{-\tau_0}),\]
which gives  \cref{CS1}.

Now, in order to prove \cref{CS2},  by restricting $U$, we can assume that $U$ is simplectomorphic to
\[ \A:=\T^{n}\times [0,1]^{n}\]
 equipped with the standard symplectic form.
 We endow $\A$ with the metric $d_\infty  ((z_\ell)_\ell ,(z_\ell')_\ell)= \max_{1\le k\le 2n} |z_\ell-  z_\ell'|$. If $d$ is the Riemannian metric of $M$, then there exists $c>0$ such that for every $z$ and $z'$ in $U$ with coordinates $(z_\ell)_\ell\in \A$ and $(z_\ell')_\ell \in \A$
 \[ d(z,z')\ge c d_\infty  ((z_\ell)_\ell ,(z_\ell')_\ell).\]
Obviously, to prove \cref{CS2} it suffices to show
\begin{lemma}\label{CS3}
There exists $\lambda>0$ such that for any $\tau<2n$ and $\epsilon_0>0$, there exist $\epsilon<\epsilon_0$ and a nonempty open set $\mathcal V_{\tau,\epsilon}\subset\Symp^\infty(\A\! \hookrightarrow\!  \T^n\times \R^n )$
near the canonical inclusion $\A \hookrightarrow \T^n\times \R^n $ such that every $g\in \mathcal V_{\tau,\epsilon}$ 
    leaves invariant an open subset $A\subset \A$ of measure $\ge \lambda$, such that the emergence at scale $\epsilon/c$ of $g|_A$ is at least $\exp(\epsilon^{-\tau})$ for $A$ endowed with the normalized measure $\leb|_A$.
\end{lemma}

\begin{proof}[Proof of \cref{CS3}] 
We shall build on the argument which was developed in \cite[Lemma 5.8]{BB21} for the two-dimensional case.  The idea is to select a number 
$\ge \exp (\epsilon^{-\tau})$ of disjoint strips in $\A$ and send them $\epsilon$  apart,  in mean, via a symplectomorphism $H$. Then we take an integrable twist map of $\A$ which leaves the strips invariant and consider the difeomorphism which is conjugate by $H$ to this twist map. By construction, the orbits in the image $H(S)$ of a strip will be well-equidistributed in $H(S)$, and so any pair of orbits that lie in the images of two different strips will be $\epsilon$-distant in mean. This will imply the sought bound on the emergence for the diffeomorphism under consideration. Furthermore, the KAM theorem is applicable here (because of the conjugacy to a non-degenerate twist), which implies the persistence of most orbits in the images of these strips, so the bound on the emergence at the scale $O(\epsilon)$ remains valid for an open set of perturbations.

This argument worked well in dimension 2 in in \cite{BB21}, however it faces a new difficulty in higher dimensions: it is not immediately clear that the sought $H$ exists (due to symplectic phenomena, such as the non-squeezing theorem). Yet, we can find a map $H$ which satisfies the desired property modulo a small set, using Katok ``Basic Lemma'' \cite{Ka73}.
\medskip 

\noindent\textbf{First step: Construction of the symplectic conjugacy $H$.}  

Let $k$ be a large even integer.  We first define a disjoint family of ``target  boxes'' $(C_{j})_{ j\in \sJ}$, where $\sJ=  
 \{1,\dots, k\}^{2n}$, such that
for every $ j = (j_i)_{1\le i\le 2n} \in  \sJ$
  \[C_{ j}:=\prod_{i=1}^{2n}\left  [\frac{j_i - \tfrac23}{k}, \frac{j_i-\tfrac13}{k} \right ]\subset \A.\]
Next, we introduce the notation
\[  M := \left[\sqrt[n] {(2\#\sJ)^{-1/4}\exp(\#\sJ/20)}\right]= \left[ (2k^2)^{-1/4}\exp\left(\frac{k^{2n}}{20n}\right)\right ], 
\]
\[
\sQ :=  \{ 1, \dots,  k^2  \}^{n-1} \times  \{ 1, \dots, \tfrac12 k^2  \} \qand  \sP:=  \{ 1, \cdots , M \} ^n.\]
Note that
 \[\# \sQ= \tfrac12 k^{2n}= \tfrac12 \#\sJ \qand \#\sP= M^{n}.\] 

Define $Y_p:=\prod_{i=1}^{n}\left  [\frac{p_i-\tfrac{15}{28} }{M}, \frac{p_i-\tfrac{13}{28}}{M} \right ]$ for $p\in\sP$, and consider the annuli $A_p=\T^n\times Y_p$. We will fix a small  $\eta>0$ in the sequel, see (\ref{etadef}). For every $p=(p_i)_i \in \sP$, let us
consider a system of disjoint boxes 
\[B_{q,p}=   \prod_{i=1}^{n-1}\left  [\frac{q_i -1+\eta }{k^2}, \frac{q_i -\eta }{k^2} \right ]\times \left  [\frac{2q_n - 2 +\eta }{k^2}, \frac{2q_n - \eta }{k^2} \right ]
\times Y_p,\]
whose union over all $q= (q_i)_i\in \sQ$ ``almost covers'' the annulus 
$A_p$ (with the accuracy $O(\eta)$).
 
We are going to construct a symplectomorphism $H$ which sends each small box $B_{q,p}$ inside some target box $C_j$ such that  the families of the images $(H(B_{q,p}))_{q\in \sQ}$ corresponding to fixed $p\in \sP$, stay uniformly apart from each other (so the images of the annuli 
$A_p$ stay uniformly apart from each other). The possibility of such arrangement is based on the following combinatorial lemma.
\begin{lemma}[{\cite{BB21}, \cite[Lem. 4.5]{Ber22}}]
\label{coloriage}For every large $k$, there is  a map
\[\Phi : \sQ\times \sP \to \sJ\] such that:
\begin{enumerate}[(1)]
\item  $\# \Phi^{-1}(\{j\}) \le    \tfrac34 \#\sP$ for every $j\in  \sJ
$;  
\item  $\#\{(q,q')\in \sQ^2: \Phi(q,p)=\Phi(q',p')\}<\tfrac34  \#\sQ$
for every pair of indices $p\neq p'$.\end{enumerate}
\end{lemma}

Let us take such map $\Phi$. For every $j\in \sJ$, the volume of
$\bigsqcup_{(q,p)\in \Phi^{-1}(\{j\})} B_{q,p}$ is at most
 \[\tfrac34 \#\sP \cdot \frac2{k ^{2 n}}\cdot \frac1{(14M)^n }=
\tfrac32  \cdot \frac1{14 ^n\cdot  k ^{2 n}} 
<\leb\,  C_j= \frac{1}{ (3 k) ^{2n}}\]
(we use here the first property of the map $\Phi$ of the lemma above). 
Hence, one can choose a system of disjoint compact sets 
$\hat B_{q,p}\subset C_{\Phi(q,p)}$ such that the volume of
$B_{q,p}$ equals to the volume of $\hat B_{q,p}$ for each 
$(q,p)\in\sQ\times\sP$,

By Katok lemma \cite[Basic lem. 3.1]{Ka73}, for every $\eta>0$ there exists $H\in \Symp^\infty(\A)$ such that 
$$\leb(\hat B_{q,p})\Delta H(B_{q,p}))<\eta\cdot  \leb(B_{q,p}),$$
and $H$ is equal to the identity near the boundary of $\A$. 

\medskip 

\noindent\textbf{Second step: A near-identity map $g$ which is $H$-conjugate to an integrable twist.}  

Let $V:=  [\tfrac1{2M}, 1-\tfrac1{2M}]^n$ and let $\rho: \A \to \R$ be a $C^\infty$ function supported by  
$int\,  \A$ and  equal to $1$ on  $\T^n\times V$. 
Let $(G^t)_{t\in \R}$ be  the flow in $\A$ defined by the Halmitonian
\[h(\theta , y ) =  \rho(\theta, y)   y^2,\]
where $\theta\in \T^n$, $y\in [0,1]^n$. The corresponding differential equation, when $(\theta, y )\in \T^n\times V$, is 
\[\dot \theta = y, \qquad \dot y =0,\]
so $G^t$ leaves $\T^n\times V$ invariant and acts as
\[(\theta,y)\mapsto (\theta + t y, y).\]
Obviously, every torus $\T_y:=\T^n\times \{y\}$ is invariant with 
respect to $G^t$ if $y\in V$. Moreover, for every $t>0$, for a.e. $y\in V$,  the corresponding frequency vector $\alpha=ty$ is Diophantine, i.e., $\T_y$ is a KAM-torus for $G^t$. Since the matrix $\partial_y \alpha = t \cdot id$ is invertible, all theses KAM-tori are non-degenerate.

Note that $G^t|_{\T_y}$ is ergodic and every its orbit is equidistributed in  $\T_y$ for any KAM-torus $\T_y$.  Thus, for every $p=(p_i)_i\in \sP$,  the $G^t$-orbit of almost every  point in the annulus $A_p$ spends an average time $\ge 2\frac{(1-2\eta)^n}{ k^{2n}}$  in the box $B_{q,p}$ for each $q\in \sQ$.  

Let $B_{q,p}'$ be the subset of $B_{q,p}$ which is sent into the target box $C_{\Phi(q,p)}$ by the symplectomorphism $H$. By the construction 
of $H$, we have $\leb\, B_{q,p}'  \ge (1-\eta)\cdot \leb\, B_{q,p}$. Hence 
there is a subset  of 
$Y_p=\prod_{i=1}^{n}
\left  [\frac{p_i-\tfrac{15}{28}}{M}, \frac{p_i-\tfrac{13}{28}}{M} \right ]$  of  measure $\ge (1-\sqrt \eta)\cdot \leb Y_p$ formed by the points $y$ such that the intersection of $B_{q,p}'$ with the torus $\T_y= \T^n \times \{y\}$  has relative measure 
$\ge (1-\sqrt \eta)$ in $B_{q,p}\cap \T_y$. Consequently, there exists a subset $Y_p'\subset Y_p$ of  measure $\ge (1-\# \sQ\cdot  \sqrt \eta) \leb Y_p$ such that, for every point in $\T^n\times  Y_p'$, its orbit by  $G^t$ spends an average time $\ge 2\frac{(1-2\eta)^n }{ k^{2n}}(1-\sqrt \eta)$ in each of the sets $B_{q,p}'$, $q\in \sQ$.

Since $G^t$ tends to identity as $t\to 0$, the map $g= H\circ G^t\circ H^{-1}$ gets as close to identity as we want if $t$ is small enough. Up to the conjugacy, the dynamics of $g$ is the same as the dynamics of $G^t$. In particular, for every $p\in \sP$, for every point in 
$H(\T^n\times Y_p')$, its orbit by $g$ spends an average time $\ge2 \frac{(1-2\eta)^n}{ k^{2n}} (1-\sqrt \eta)$  in each of the sets $H(B_{q,p}')\subset C_{\Phi(q,p)}$, $q\in \sQ$. 

\medskip 

\noindent\textbf{Third step: Estimate for the emergence of the $H$-conjugate twist map $g$.}  

Let us now fix $\eta$ depending on $k$ such that 
\begin{equation}\label{etadef}
  {(1-2\eta)^n}   (1-\sqrt \eta)>1-\frac{1}{60k},\qquad (1-\#\sQ\cdot\sqrt  \eta)>1-\frac{1}{60k}.
\end{equation}
 By the construction of $g$,
for each annulus $A_p$, $p\in\sP$, there is a subset $A'_p\subset H(A_p)$ of measure  $(1-\frac{1}{60k})\cdot \leb\, A_{p}$  which consists of points $z$ such that their corresponding empirical measure with respect to the map $g$  can be represented as
\[e(z) = \frac{1}{60k} \cdot \mu_0^z + \frac{1-\tfrac{1}{60k}}{\# \sQ} 
  \sum_{q\in \sQ}   \mu_{\Phi(q,p)} ^z,\]
where $\mu_{\Phi(q,p)} ^z$ is a probability measure supported by the target box $C_{\Phi(q,p)}$, and $\mu_0^z$ is a probability measure on $\A$.

For every $  j\neq j' \in \sJ$, we notice that the centers $z^j$ and $z^{j'}$  of the boxes $C_j$ and $C_{j'}$ are at least $\frac1k $-distant.
For every $p\in \sP$, let $\mu_p$ be the probability measure equidistributed on $(z^{\Phi(q,p)})_{q\in \sQ}$:
\[\mu_p:= \frac 1{\# \sQ}  \sum_{j= \Phi(q,p); q\in \sQ}   \delta_{z^{j} }.\]  
Since the box  $C_j$ is the $\tfrac1{6k}$-ball about $z_j$, the measure $\mu_{\Phi(q,p)}^z $ is $\tfrac1{6k}$ close to   
$\delta_{z^{j}}$, for every $z\in A_p$. Hence, for all $z\in A'_p$,
\[d(\mathsf e(z), \mu_p) \leq \frac { 1-\tfrac{1}{60k} }{6 \cdot k} +\frac{1}{60k} \cdot   \diam \A < \frac{11}{60 \cdot k}\]
(recall that $\diam A=1$). It follows that for every $p\in \sP$,
\[ \frac1{\leb(\A_p)} \int_{H(\A_p)}  d(\mathsf e(z), \mu_p)d\leb <
 \left(1-\frac{1}{60k}\right)  
\frac { 11}{60 \cdot k} + \frac{1}{60\cdot k}  < \frac {1}{5 \cdot k}.\]
  
Thus, the normalized ergodic decomposition of  $\mathsf e|_{H(\bigsqcup_{p\in \sP} A_p)}$ is $\frac{1}{5 k} $-close to
\[\hat \mu:= \frac1{\#\sP}  \sum_{p\in \sP}  \delta_{\mu_p}.\] 
\begin{claim}\label{cl.apart}For any $p\neq p' \in \sP$
\[d(\mu_p, \mu_{p'})\ge \frac 1{4k}.\] 
\end{claim}
\begin{proof}[Proof of the claim]
For $1\le j\neq j'\le N$, the distance between $z^{j}$ and $z^{j'}$ is at least $1/ {k}$. By the second property of the map $\Phi$ in \cref{coloriage}, for any transport plan from $\mu_p$ to $\mu_{p'}$ there has to be at least $\# \sQ/4$ masses of weight $1/\#\sQ$ which will be moved
from a certain $z^{j}$ different to $z^{j'}$. So the transport cost of each will be at least $1/(k\#\sQ)$. The summation implies  the bound.
\end{proof}

As we see, the probability measure $ \hat \mu$ on the space of probability measures on $\A$ is equidistributed on $\#\sP$ atoms which are at least $\frac 1{4k}$-distant. Therefore, the distance from $\hat \mu$ to any measure supported by $\#\sP/10$ atoms is larger than 
$\frac {0.9}{4\cdot  k}$, see \cite[Lemma 3.19]{BB21}. 
Thus, the distance between the normalized ergodic decomposition of $\mathsf e|_{H(\bigsqcup_{p\in \sP}  A_p)}$ to a measure supported by $\#\sP /10$ atoms is at least  $\frac {0.9}{4\cdot k}- \frac 1{5\cdot k}= \frac 1{40\cdot k}$. 

This means that the emergence of $g|_{H(\bigsqcup_{p\in \sP}  A_p)}$
 at scale $\frac 1{40\cdot k}$ is at least $\# \sP/10$. 
Observe:
 \[  \frac{\# \sP}{10 } = \frac{M^n}{10 }\sim \tfrac1{10} (2k^2)^{-n/4} \cdot \exp (k^{2n} /20) \gg \cdot \exp (  (c\cdot  \tfrac{1}{40k})^{-\tau}),\]
as $\tau <2n$ and $k$ is taken very large. 
So, the emergence at scale $\epsilon/c$ for the map $g|_{H(\bigsqcup_{p\in \sP}  A_p)}$ is greater than $2\cdot\exp( \epsilon^{-\tau})$ for $\epsilon= \frac{c}{40k}$.\medskip

\noindent\textbf{Conclusion.}
Recall that the set $H(\bigsqcup_{p\in \sP}  A_p)$ is foliated by invariant tori of $g$, and a.e torus is a non-degenerate KAM-torus.
By the Lazutkin-P\"{o}schel version of the KAM theorem \cite{Laz73,Po82}, 
``the most'' of these tori persist for all exact symplectomorphisms 
which are $C^\infty$-close to $g$. 

Namely,  there exists  an  invariant set $A\subset H(\bigsqcup_{p\in \sP}  A_p)$ of relative measure close to $1$, formed by the KAM-tori  which persists for any  exact symplectic $C^\infty$-perturbation of $g$: the persistent   KAM tori  
depends on the perturbation continuously. 
 Thus, the emergence at scale $\epsilon/c$ gets close to the emergence at scale $\epsilon/c$ for the restriction to $H(\bigsqcup_{p\in \sP}  A_p)$, i.e., it is at least $\exp (\epsilon^{-\tau})$ for an open set
of maps.  \end{proof}

\section{Universality is a localizable property}\label{Universality proof}
Here, we prove \cref{prop open renorm}. Notice that it is a direct consequence of the following
\begin{proposition}\label{id+}  
For any nonempty open $\cU\subset \Symp^\infty(\B^{2n}\!\hookrightarrow\! \R^{2n})$ (where $\B^{2n}$ is the unit ball in $\R^{2n}$, in any neighborhood of the identity map in $\Symp^\infty(\B^{2n})$ there exists a  compactly supported  $\phi\in \Symp^\infty(\B^{2n}) $  such that a renormalized iteration of $\phi$ belongs to $\cU$.
\end{proposition}

\begin{proof}
Let $\xi \in C^\infty(\R, \R) $ denote a  function with  support in $(-1,1)$ such that $\xi$ is identically $1$ on $[-\tfrac12, \tfrac12]$. 
We denote by $(x,y)=(x_1,\dots, x_n,y_1, \dots, y_n) $ the  canonical coordinates on $\R^{2n}$, so the symplectic form is $\sum_k dx_k \bigwedge dy_k$.\\

\noindent{\bf  First perturbation: creation of a homoclinic connection.} 
For  any $C^\infty$-function $H$ on $\mathbb{R}^2$, when $\tau $ is small, the time-$\tau$ map $F$ by the flow generated by the Hamiltonian
\begin{equation}\label{1}
\xi( \|(x ,y )\|) \cdot \sum_k   H(x_k,y_k)
\end{equation}
is a small perturbation of the identity in $\Symp^\infty(\B^{2n})$; the perturbation is compactly supported in $\B^{2n}$.

We choose the 1-degree-of-freedom Hamiltonian $H$ such that the corresponding planar flow has a hyperbolic equilibrium at zero, with a homoclinic loop that lies entirely in the disc $\|(q,p)\|<\tfrac12 $.
For convenience of computations, we also assume that $H(q,p)=pq$ in an $\varepsilon$-neighborhood of $(q,p)=0$, i.e., the Hamiltonian vector field in the $\varepsilon$-neighborhood of $(q,p)=0$ is
$$\dot q = -q, \qquad \dot p = p.$$
The flow map in this neighborhood is linear:
$$(q,p) \mapsto (e^{-t} q, e^{t} p);$$
the interval $(q=0, |p|<\varepsilon)$ is a local unstable manifold and the interval $(|q|<\varepsilon, p=0)$ is a local stable manifold.
The existence of the homoclinic loop means that one can choose a pair of  points, $(q^s,0)$ in the local stable manifold and $(0,p^u)$ in the local unstable manifold, such that these points belong to the same orbit of the Hamiltonian flow (the homoclinic loop). Moreover, we choose these points such that the flight time between $(0,p^u)$
 and $(q^s,0)$ is an integer multiple of $\tau$, i.e.,  
 $(q^s,0)$ is the image of $(0,p^u)$ by the $m^{th}$ iteration of the time-$\tau$ map for some integer $m>0$. A piece of the local unstable manifold around $(0,p^u)$ and a piece of the local stable manifold around $(q^s,0)$ belong to the same homoclinic orbit, so the stable piece is the image of an unstable piece by the same $m^{th}$ iteration of the time-$\tau$ map.

Returning to the Hamiltonian (\ref{1}), we obtain that the time-$\tau$ map $F$ defined by it has a hyperbolic fixed point at the origin $O$ of $\R^{2n}$  with the local unstable and stable manifolds:
\[W^u_{loc}= \{ (x,y)\in \B^{2n}: x=0, \;   \|y\|<\varepsilon\},\qand W^s_{loc}=  \{ (x,y)\in \B^{2n}:  \|x\|<\varepsilon,\;  y=0\}.\]
 The map $F$
in the $\varepsilon$-neighborhood of $O$ is given by
\begin{equation}\label{4}
(x,y) \mapsto (e^{-\tau} x, e^{\tau} y).
\end{equation}
Denote $x^s=(q^s,\ldots, q^s)\in \mathbb{R}^n$ and $y^u=(p^u,\ldots, p^u)\in \mathbb{R}^n$.  The points
\[ M^s:= ( x^s, 0)\in W^s_{loc} \qand M^u:= (0,  y^u)\in W^u_{loc}\]
 belong to the same homoclinic orbit: $M^s = F^m (M^u)$. Moreover, this is an orbit of a flat tangency: the image by $F^m$ of some ball in $W^u_{loc}$ around $M^u$ lies inside $W^s_{loc}$. So, the map $F^m: (x,y)\mapsto (\bar x,\bar y)$ can be written near $M^u=(0,y^u)$ as
\begin{equation}\label{2}
 \bar x = x^s + a x + b (y-y^u) + \varphi_2(x,y-y^u) , \qquad \bar y = x\cdot  (c +  \varphi_1(x,y-y^u))\end{equation}
where $a,b,c$ are some scalars; the function $\varphi_1$ vanishes at zero, and $\varphi_2$ vanishes at zero along with its derivative. Since $F^m$ preserves the standard symplectic form, one notices that
\begin{equation}\label{6}
bc = -1.
\end{equation}

\noindent {\bf  Second perturbation: splitting the flat homoclinic tangency.} Take some small $\rho>0$ -- such that the $\rho$-neighborhood of the homoclinic point $M^s$ contains no other points of the orbit of $M^s$. Take a small $\nu>0$, and take any $C^\infty$-function $V: \mathbb{R}^n \to \mathbb{R}$.
At the point  $(x^s+x,y)$, the vector field defined by the Hamiltonian  
$$H(x^s+x,y) = \xi(\tfrac{\|(x ,y)\|}\rho) \cdot  (V( \tfrac x \nu ) - (by^u+x^s)\nu^{-2} x)$$
 is given by
$$\begin{array}{l} \displaystyle \dot x = (V(\tfrac x \nu ) - (by^u+x^s) \nu^{-2} x) \cdot \partial_y [ \xi(\tfrac{\|(x ,y)\|}\rho)], \\ \\
\displaystyle
\dot y = \xi(\tfrac{\|(x ,y)\|}\rho) \cdot ((by^u + x^s)\nu^{-2} -\nu^{-1} \nabla V(\tfrac x \nu) )  + (V(\tfrac x \nu )  - (by^u+x^s) \nu^{-2} x)\cdot \partial_x [\xi(\tfrac{\|(x ,y)\|}\rho)].\end{array}$$
It vanishes outside the $\rho$-neighborhood of $M^s$ and its $C^r$-norm is dominated by $ \nu^{-(r+1)} $ for any $r\ge 1$. Hence, the corresponding time-$\nu^{r+2}$ map $\phi_{r,\nu}$ is identity outside the $\rho$-neighborhood of $M^s$ and is $O(\nu)$-close to identity in $C^r$. In the radius $\rho/2$ ball around $M^s$ the vector field at $(x^s+x, y)$ is given by
$$\dot x = 0, \qquad
\displaystyle
\dot y = (by^u + x^s)\nu^{-2} -\nu^{-1} \nabla V(\tfrac x \nu ), $$
so in a somewhat smaller neighborhood of $M^s $ (which does not shrink as $\nu\to 0$) the map  $\phi_{r,\nu}$ is given by
$$(x^s+  x,y) \mapsto ( x^s+x, \; y + (by^u + x^s)\nu^r - \nu^{r+1} \nabla V(\tfrac x \nu)).$$

Thus, if $\nu$ is small enough, the symplectic map $\tilde F=\phi_{r,\nu}\circ F$ is a $C^r$-small perturbation of $F$, so it is a small perturbation of the identity. Since the perturbation is localized in the $\rho$-neighborhood of $M^s$, we have $\tilde
 F^m=\phi_{r,\nu} \circ F^m$ in a neighborhood of $M^u$, so formula (\ref{2}) is perturbed to
\begin{equation}\label{3}\begin{array}{l}\displaystyle
\bar x = x^s + a x + b (y-y^u) + \varphi_2(x, y-y^u),\\ \\ \displaystyle
\bar y =  x\cdot (c +  \varphi_1(x,y-y^u) )+(by^u + x^s)\nu^r - \nu^{r+1} \nabla V(\nu^{-1} (\bar x - x^s)). 
\end{array}
\end{equation}

\noindent {\bf Homoclinic renormalization.} We choose the small $\nu$ such that there exists a large integer $j$ satisfying
\begin{equation}\label{7}
b \nu^r=e^{-j \tau}.
\end{equation}
Consider the box near $M^s$ equal to:
$$B_j=\{\|x-x^s\|< 2\nu, \|y-e^{-j\tau} y^u \| < 2\nu^{r+1} \}.$$
Let $X=x-x^s, Y=e^{j\tau}y -y^u$. In the new coordinates, the box is given by
$$B_j=\{\|X\|< 2\nu, \|Y\| < \frac{2}{b} \nu\}.$$
By (\ref{4}), 
 the image $F^j B_j$ lies in an $O(\nu)$-small neighborhood of the homoclinic point $M^u$: 
the point $(X,Y)$ is sent by $F^j$  to
\[(x,y) =   (e^{-j\tau}(x^s+ X),   y^u +Y)  
.\]
Hence, the image by $\tilde F^{j+m}$ is given, according
to (\ref{3}), by
\[\begin{array}{l}\displaystyle
\bar x = x^s + a e^{-j\tau}(x^s+ X) + b  Y + \varphi_2(e^{-j\tau}(x^s+ X), Y),\\ \\ \displaystyle
\bar y = e^{-j\tau}(x^s+ X)  (c  +  \varphi_1(e^{-j\tau}(x^s+ X),Y) )+ (by^u + x^s)\nu^r - \nu^{r+1} \nabla V(\nu^{-1} (\bar x - x^s)).
\end{array}
\]
We have $\varphi_i(e^{-j\tau}(x^s+ X), Y) = \varphi_i(0, Y)+O(e^{-j\tau})$ (in the $C^r$-norm), so 
\[\begin{array}{l}\displaystyle
\bar x = x^s    + b  Y +    \varphi_2(0, Y) +O( e^{-j\tau}),\\ \\ \displaystyle
\bar y =  e^{-j\tau}(x^s+ X)  (c  +  \varphi_1(0,Y) ) + e^{-j\tau}(
y^u + x^s/b)- e^{-j\tau} \frac{\nu}b \nabla V(\nu^{-1} (\bar x - x^s)) + O( e^{-2j\tau}).
\end{array}\]\\
 
Thus, the map $\tilde F^{j+m}$ takes $B_j$ to a small neighborhood of $M^s$ and has the form $(X,Y)\mapsto (\bar X, \bar Y)$ where
\begin{equation}\label{5}\begin{array}{l}\displaystyle
\bar X = b Y +    \varphi_2(0, Y)+  O(\nu^r ) ,\\ \\ \displaystyle
 \bar Y =   c  X + x^s(c+1/b)+(x^s+ X)   \varphi_1(0,Y) - \tfrac{\nu}b   \nabla V(\nu^{-1} \bar X) +
 O(\nu^r ). \end{array}
\end{equation}
By (\ref{6}) we have $(c+1/b)=0$. Also, since $\phi_1(0)=0$, $\phi_2(0)=0$, and $D\phi_2(0)=0$, one sees that after the scaling $(X,Y)\to (\nu X, \nu Y/b)$ the map (\ref{5}) becomes 
\[\bar X =  Y + O(\nu), \qquad  
\bar Y = -X + C Y - \nabla V(Y)    +  O(\nu),\]
with the definition domain $\max\{\|X\|,\|Y\|\}\leq 2$.
We denoted here $C={x^s}\cdot  \partial_y \phi_1(0)$;
the $O(\nu)$ terms here are $O(\nu)$ in $C^r$. 

Note that
the coordinate transformation we made were affine and conformally symplectic, so the map preserves the standard symplectic form $dX\wedge dY$ for all $\nu>0$.  Thus, 
$$dX\wedge dY = d\bar X\wedge d\bar Y=  - dY\wedge dX
+ dY \wedge C dY - dY \wedge (\partial_{YY} V) dY + O(\nu).$$
Since the Hessian matrix $\partial_{YY} V$ is symmetric, we have $dY \wedge (\partial_{YY} V) dY=0$, so taking
the limit $\nu\to 0$, gives that $dY \wedge C dY=0$, i.e., the matrix $C$ must be symmetric, so $C Y = \frac12 \nabla <\!Y,CY\!>$.

Thus, we have shown that the perturbation $\tilde F$ can be chosen such that some renormalized iteration of it gets as close as we want to the Henon-like map
\[\bar X =  Y, \qquad \bar Y = - \nabla \hat V(Y)  -  X,\]
for any given function $\hat V$ (here, $\hat V = V - \frac12  <\!Y,CY\!>$). By \cite{Tu03}, the union of the sets of renormalized iterations of Henon-like maps with all possible $C^\infty$ functions $\hat V$ is dense in $\Symp^\infty(\B^{2n}\!\hookrightarrow\!\R^{2n})$. Since a renormalized iteration of a renormalized iteration of $\tilde F$ is again a renormalized iteration of $\tilde F$, we obtain that given any $g \in \Symp^\infty(\B^{2n}\!\hookrightarrow\! \R^{2n})$, one can choose $\tilde F$ such that some of its renormalized iterations would approximate $g$ with any given accuracy in $C^\infty$.
\end{proof}

\end{document}